\magnification=\magstep1
\input amstex
\documentstyle{amsppt}
\hoffset=-.35truein
\hsize=7.25truein
\vsize=8.95truein

\topmatter
\title
A Hilbert-Schmidt integral operator and the Weil distribution
\endtitle
\abstract
 In this paper, a positive operator is given. It is shown that the product
of this positive operator and the convolution operator is a trace class
Hilbert-Schmidt integral operator and has nonnegative eigenvalues.
A formula is given for the trace of this product operator.
It seems that this product operator is the closest trace class
integral operator which has nonnegative eigenvalues and
is related to the Weil distribution.
A relation is given between the trace of the
product operator and the Weil distribution.
\endabstract
\author
Xian-Jin Li
\endauthor
\address
Department of Mathematics, Brigham Young University, Provo, Utah 84602, USA
\endaddress
\email
xianjin\@mathematics.byu.edu
\endemail
\endtopmatter
\document

\vskip 0.15in
\heading
1. Introduction
\endheading

    We first introduce some notations.
  Let $C_c^\infty(0, \infty)$ be the space of compactly supported
  smooth functions on $(0,\infty)$. Because of [1, Theorem 1, p. 191],
  we assume that $g$ is a real-valued function in $C_c^\infty(0, \infty)$.
   Define
 $$h(x)=\int_0^\infty g(xy)g(y)dy.\tag 1.1$$
  Then $h\in C_c^\infty(0, \infty)$.
   Since $g$ belongs to $C_c^\infty(0, \infty)$, there exists a positive
   number $\mu<1$ such that the support of $g$
is contained in $(\sqrt\mu, \mu^{-1/2})$.  Thus,
$h(x)=0$  for all $x\not\in (\mu, \mu^{-1})$.

    Let $\Lambda>1$ be any fixed number, and let
\roster
\item $S=\{\infty$, all primes $p<\mu^{-1}\}$, $S^\prime=S-\{\infty\}$,
$Q=$ field of rational numbers, $Q^*=Q-\{0\}$,

\item $\Bbb  A_S=\Bbb R\times\prod_{p\in S^\prime} Q_p$,
$J_S=\Bbb R\times \prod_{p\in S^\prime}Q_p^*$, $C_S=J_S/{O_S^*}$,
$Q_p=$ $p$-adic completion of $Q$,

\item $O_S^*=\{\xi\in Q^*:\,\, |\xi|_p=1,\,\, p\not\in S\}$,
 $I_S=\Bbb R_+\times\prod_{p\in S^\prime}O_p^*$,
 $O_p^*=\{x\in Q_p:|x|_p=1\}$,

\item $E_S(f)(x)=\sqrt{|x|}\sum_{\xi\in O_S^*}f(\xi x)$,
 $e_S(f)(x)=\sum_{\xi\in O_S^*}f(\xi x)$,

\item $S_\Lambda$ is the multiplication operator multiplied by the function
which equals $0$ when $|x|\leq \Lambda^{-1}$ and equals $1$
when $|x|>\Lambda^{-1}$,

\item $P_\Lambda$ is the multiplication by the function which is $0$ for
$|x|\geq\Lambda$ and is $1$ for $|x|<\Lambda$,

\item $\Bbb N_S$ is the set of all positive integers which are
 products of primes in $S^\prime$,

 \item $h(x):= h(|x|)$  for $x\in\Bbb A_S$ or $C_S$, $|x|:=|x|_S$.
\endroster

Let $\psi_p$ be the additive character on $Q_p$ given as in [10, Theorem 2.2.1, p. 309].
  Let $\Psi_S=\prod_{p\in S}\psi_p$.  For $f=\prod_{p\in S} f_p\in L^2(\Bbb  A_S)$,
we define
$$\frak F_Sf(\beta)=\int_{\Bbb A_S}f(\alpha)\Psi_S(-\alpha\beta)d\alpha.$$
Then $\frak F_Sf=\prod_{p\in S}\frak F_pf_p$.  As the Schwartz-Bruhat space
$S(\Bbb A_S)$ is dense in $L^2(\Bbb A_S)$,
the definition of $\frak F_S$ can be extended to all functions in $L^2(\Bbb A_S)$.
The set of all functions $f$ in $S(\Bbb A_S)$ satisfying $f(0)=0$
and $\frak F_Sf(0)=0$ is denoted by $S(\Bbb A_S)_0$.

The convolution operator $V(h)$ is defined by
$$V(h)F(x)=\int_{C_S}h(\lambda^{-1}x)\sqrt{|\lambda^{-1}x|}
F(\lambda)d^\times\lambda $$
for $F\in L^2(C_S)$.
  The Weil distribution $\Delta(h)$ [13, p. 18] is given by
  $$\Delta(h)=\widehat h(0)
+\widehat h(1)-\sum_p\int_{Q_p^*}^\prime
{h(|u|_p^{-1})\over |1-u|_p} d^*u, $$
where the sum on $p$ is over all primes of $Q$ including the infinity
prime, and for $p\neq\infty$
$$\int_{Q_p^*}^\prime{h(|u|_p^{-1})\over |1-u|_p} d^*u
=\sum_{m=1}^\infty \log p\left[h(p^m)+p^{-m}h(p^{-m})\right].$$
In particular,
 $$\int_{Q_p^*}^\prime {h(|u|_p^{-1})\over |1-u|_p}d^*u=0 $$
 for $p\not\in S$ by our choice of $S$.

  Similarly as function field case [6] we always assume that
  $\widehat h(1)=\widehat h(0)=0$; see [5, Theorem 1.3].

  The paper is organized as follows:  Let
    $$T=S_\Lambda\left(S_\Lambda
    - E_S\frak F_S^t P_\Lambda\frak F_S E_S^{-1}\right)S_\Lambda.$$
   We show in Section 2 that $T$ is a positive operator on $L^2(C_S)$
   and that all eigenvalues of $V(h)T$ are nonnegative, namely the following
   two theorems.

  \proclaim{Theorem 1.1}  Let $Z_\Lambda =\frak F_S^t P_\Lambda\frak F_S$ and
$$T=S_\Lambda\left(S_\Lambda- E_SZ_\Lambda E_S^{-1}\right)S_\Lambda.$$
Then $T$ is a positive operator on $L^2(C_S)$.
\endproclaim

\proclaim{Theorem 1.2} All eigenvalues of $V(h)T$ on $L^2(C_S)$ are nonnegative.
\endproclaim

   In Section 3, we prove that $V(h)T$ is a trace class operator on $L^2(C_S)$.

\proclaim{Theorem 1.3} $V(h)T$ is of trace class on $L^2(C_S)$.
\endproclaim

 The trace of the operator $V(h)T$ is computed in Section 4.
 In particular, we obtained the following theorem.

  \proclaim{Theorem 1.4}   We have
$$\aligned \text{trace}_{L^2(C_S)}\left\{V(h)T\right\}
=&\int_{\Bbb A_S, |v|>1}\Psi_S(v)
\log|v|dv\int_{\Bbb A_S, |u|\geq 1}h(u)\Psi_S(-uv)du\\
&+\int_{|v|>1}\Psi_S(v)\log|v|dv\int_{|u|>1}h(u)\Psi_S(-uv)du.
\endaligned$$
\endproclaim

By Theorem 1.2, all eigenvalues of $V(h)T$ on $L^2(C_S)$ are nonnegative.
Theorem 1.3 says that $V(h)T$ is of trace class on $L^2(C_S)$.
The Lidski\u{\i} theorem (see Lemma 5.2) asserts that spectral and functional
traces of a trace class operator are the same.  So, the trace of $V(h)T$
is nonnegative.   The following corollary then follows from Theorem 1.4;
see Section 5.

 \proclaim{Corollary 1.5}  Let $h$ be any function given as in (1.1).  Then
 $$\sum_{k_1, k_2, l_1, l_2\in\Bbb N_S}{\mu(k_1)\mu(k_2)\over k_1k_2}
  \int_1^\infty \log v\,\cos{2\pi l_1v\over k_1}dv
  \int_1^\infty h(u)\cos{2\pi l_2uv\over k_2}du\geq 0.$$
 \endproclaim

 We derive in Section 6 the following relation between the trace of
 the operator $V(h)T$ and the Weil distribution $\Delta(h)$.

 \proclaim{Theorem 1.6} We have
$$\aligned \text{trace}_{L^2(C_S)}\left\{V(h)T\right\}
=\Delta(h)&-\int_{|v|<1}\log|v|\Psi_S(-v)dv\int_{|u|\geq 1}h(u)\Psi_S(uv)du\\
&-\int_{|v|<1}\log|v|\Psi_S(-v)dv\int_{|u|>1}h(u)\Psi_S(uv)du.
\endaligned$$
\endproclaim

\vskip 0.35in
\heading
2.  Nonnegativity of eigenvalues of $V(h)T$
\endheading

  In this section we show that all eigenvalues of $V(h)T$ are nonnegative.

   \proclaim{Lemma 2.1} ([5, Lemma 3.3 and Lemma 3.4, p. 2467])  $E_S\frak F_S E_S^{-1}$
  is unitary on $L^2(C_S)$, and its inverse is $E_S\frak F_S^t E_S^{-1}$.
\endproclaim

\demo{Proof of Theorem 1.1}  Since $S_\Lambda S_\Lambda=S_\Lambda$, we have
$$T=S_\Lambda\left(1- E_SZ_\Lambda E_S^{-1}\right)S_\Lambda.$$
By Lemma 2.1, $E_S\frak F_S E_S^{-1}$
is a unitary operator on $L^2(C_S)$.  This implies that $E_SZ_\Lambda E_S^{-1} $ is
 an orthogonal projection on $L^2(C_S)$ as it is a nonzero self-adjoint idempotent.
 So is $1- E_SZ_\Lambda E_S^{-1} $.  Thus,
$$\langle (1- E_SZ_\Lambda E_S^{-1} )F,F\rangle\geq 0$$
for $F\in L^2(C_S)$.  As $S_\Lambda$ is bounded
linear operators on $L^2(C_S)$, we have
$$\langle (1- E_SZ_\Lambda E_S^{-1} )S_\Lambda F,S_\Lambda F\rangle\geq 0$$
for $F\in L^2(C_S)$. It follows that
$$\aligned \langle TF, F\rangle
&=\langle S_\Lambda\left(1- E_SZ_\Lambda E_S^{-1} \right)S_\Lambda F, F\rangle\\
&=\langle (1- E_SZ_\Lambda E_S^{-1} )S_\Lambda F,S_\Lambda F\rangle\geq 0
\endaligned $$
for every element $F\in L^2(C_S)$.  Hence, $T$ is a
positive operator on $L^2(C_S)$.

   This completes the proof of Theorem 1.1.
\qed\enddemo

  \proclaim{Lemma 2.2} $V(h)$ is a positive operator on $L^2(C_S)$.
    \endproclaim

  \demo{Proof}   Assume that $F$ is any element in $L^2(C_S)$
  with compact support.  By definition,
   $$ V(h) F(x)=\int_0^\infty  F(\lambda)\sqrt{|x/\lambda|}d^\times\lambda
 \int_0^\infty g(|x/\lambda|y)g(y)dy.$$
 If we make the change of variables $y\to |\lambda|y$ then
 $$\aligned &\int_{C_S} V(h) F(x) \bar F(x)d^\times x\\
 &=\int_{C_S} \bar  F(x)d^\times x
 \int_{C_S}  F(\lambda)\sqrt{|x/\lambda|}d^\times\lambda
 \int_0^\infty g(|x/\lambda|y)g(y)dy\\
 &=\int_{C_S} \bar  F(x)\sqrt{|x|}d^\times x
 \int_{C_S}  F(\lambda)\sqrt{|\lambda|}d^\times\lambda
 \int_0^\infty g(|x|y)g(|\lambda| y)dy.
 \endaligned $$
  Since the above triple integral is absolute integrable as $F$ and $g$
  are compactly supported,   we can change order of integration to derive
   $$\int_{C_S} V(h)F(x)\bar F(x)d^\times x
   =\int_0^\infty \overline{(\int_{C_S}
   F(x)g(|x|y)\sqrt{|x|}d^\times x)}(\int_{C_S} F(\lambda)
  g(|\lambda| y)\sqrt{|\lambda|}d^\times\lambda)dy\geq 0. $$
  Since $S$ is a finite set, compactly supported functions in
  $L^2(C_S)$ are dense in $L^2(C_S)$.
  As $V(h)$ is bounded, we have
  $$\langle V(h)F,\rangle F\geq 0$$
  for all $F\in L^2(C_S)$.  That is, $V(h)$ is a positive operator
  on $L^2(C_S)$.

    This completes the proof of the lemma. \qed\enddemo

  \demo{Proof of Theorem 1.2}  Since $T$ is a positive operator on the
complex Hilbert space $L^2(C_S)$ by Lemma 2.2, [8, Theorem VI.9 on p. 196 and 4th
paragraph on p. 195] implies that there exists a positive bounded operator $B$ on
$L^2(C_S)$, which is self-adjoint, satisfying that $T=B^2$.

Let $\lambda$ be any eigenvalue of $V(h)T$ on $L^2(C_S)$.
Then a nonzero element $F\in L^2(C_S)$ exists such that
$V(h)T(F)=\lambda F$.  If $B(F)=0$, then
$\lambda F=V(h)TF=V(h)B\cdot B(F)=0$.  Hence $\lambda=0$ as $F\neq 0$.
Thus it suffices to consider the case when $B(F)\neq 0$.

  As $B$ is bounded, we have $B(F)\in L^2(C_S)$.
  Assume that $B(F)\neq 0$.  Since $V(h)T(F)=\lambda F$,
we have $BV(h)T(F)=\lambda B(F)$.  Thus
$$\lambda \langle B(F), B(F)\rangle_{L^2(C_S)}=\langle \lambda B(F), B(F)\rangle
=\langle BV(h)T(F), B(F)\rangle.$$
Since $B$ is self-adjoint and $T=B^2$, we have
$$\langle BV(h)T(F), B(F)\rangle=\langle V(h)T(F), B^2(F)\rangle
=\langle V(h)T(F), T(F)\rangle.$$
By Lemma 2.2, $V(h)$ is a positive operator on
the space $L^2(C_S)$.  Hence
$$\langle V(h)T(F), T(F)\rangle_{L^2(C_S)}\geq 0.$$
Therefore
$$\lambda \langle B(F), B(F)\rangle_{L^2(C_S)}\geq 0.$$
Since $B(F)\neq 0$, we must have $\langle B(F), B(F)\rangle_{L^2(C_S)}>0$
as $L^2(C_S)$ is a Hilbert space. It follows that $\lambda\geq 0$.

 This completes the proof of Theorem 1.2. \qed\enddemo

\vskip 0.35in
\heading
3.  Traceability of $V(h)T$
\endheading

  In this section we show that $V(h)T$ is of trace class on $L^2(C_S)$.

\proclaim{Lemma 3.1} ([5, (3.6), p. 2470]) We have
$$V(g)V(g^*)=V(h)=V(g^*)V(g)$$
on $L^2(C_S)$, where $g^*(x)=x^{-1}g(x^{-1})$.
\endproclaim

  We denote $[A,B]=AB-BA$.

\proclaim{Lemma 3.2}  Let $T$ be given as in Theorem 1.1.  Then
$$V(h)T=-V(g)[S_\Lambda, V(g^*)]
 S_\Lambda+V(g)[S_\Lambda,
 V(g^*)S_\Lambda E_S \frak F_S^tJE_S^{-1}]
E_S J\frak F_S E_S^{-1} S_\Lambda.$$
where $Jf(x)=|x|^{-1}f(x^{-1})$.
\endproclaim

\demo{Proof} Since $JS_\Lambda J=P_\Lambda$ and
$$ T=S_\Lambda(S_\Lambda-E_S  \frak F_S^t
P_\Lambda\frak F_S E_S^{-1})S_\Lambda,$$
 we have
$$ V(h)T=V(h)S_\Lambda(S_\Lambda-E_S  \frak F_S^t
 JS_\Lambda J\frak F_S E_S^{-1}) S_\Lambda.$$

 By Lemma 2.1, $E_S\frak F_S^t\frak F_SE_S^{-1}$ is the identity map
 on $L^2(C_S)$. Since $E_S^{-1}S_\Lambda E_S=S_\Lambda$,
 from Lemma 3.1 we derive that
$$\aligned &V(h)S_\Lambda(S_\Lambda-E_S  \frak F_S^t
  JS_\Lambda J\frak F_S E_S^{-1}) S_\Lambda\\
 &=-V(g)[S_\Lambda, V(g^*)]S_\Lambda+V(g)[S_\Lambda,
 V(g^*)S_\Lambda E_S \frak F_S^tJE_S^{-1}]
E_S J\frak F_S E_S^{-1}S_\Lambda.\endaligned $$
Hence,
$$V(h)T=-V(g)[S_\Lambda, V(g^*)]
 S_\Lambda+V(g)[S_\Lambda,
 V(g^*)S_\Lambda E_S \frak F_S^tJE_S^{-1}]
E_S J\frak F_S E_S^{-1} S_\Lambda.$$

This completes the proof of the lemma.
\qed\enddemo

 \proclaim{Lemma 3.3} ([8, Theorem VI.22(h), p. 210])
 A linear operator on a Hilbert space is of trace class if and only
 if it is a product of two Hilbert-Schmidt operators.
 \endproclaim

 \proclaim{Lemma 3.4} ([8, Theorem VI.19, p. 207])  Let
 $A,B$ be bounded linear operators on a Hilbert space $\Cal H$.
 If $A$ is of trace class on $\Cal H$, so are $AB$ and $BA$.
 If $A,B$ are of trace class, so is $A+B$.
 \endproclaim

 \proclaim{Lemma 3.5} ([8, Theorem VI.23, p. 210])  Let $\langle M, \mu\rangle$
 be a measure space.  A bounded linear operator $A$ on $L^2(M, d\mu)$ is
 Hilbert-Schmidt if and only if there exists a function
 $k(x,y)\in L^2(M\times M, d\mu\otimes d\mu)$ such that
 $$(Af)(x)=\int_M k(x,y)f(y)d\mu(y)$$
 for all $f\in L^2(M,d\mu)$.
 \endproclaim

  \proclaim{Lemma 3.6} We have
 $$ V_S(g^*)S_\Lambda E_S \frak F_S^tJE_S^{-1}F(x)
=\int_{C_S}\left(\sqrt{|x/y|}\int_{\Bbb A_S}g(\lambda)
S_\Lambda(\lambda x)\Psi_S(\lambda x/y)d\lambda\right)F(y)d^\times y$$
for $F\in L^2(C_S)$.
 \endproclaim

 \demo{Proof}   Let $F=E_S(f)$ with $f\in S(\Bbb A_S)_0$.
 We denote $J_0F(x)=F(1/x)$. Since
 $$E_S^{-1}J_0F(x)=E_S^{-1}[\sum_{\xi\in O_S^*}{1\over\sqrt{|\xi x|}}
 f({1\over\xi x})]={1\over|x|}f({1\over x})=(JE_S^{-1}F)(x),\tag 3.1$$
 by Lemma 2.1 we have
 $$\aligned &V_S(g^*)S_\Lambda E_S \frak F_S^tJE_S^{-1}F(x)
 =V_S(g^*)S_\Lambda E_S \frak F_S^tE_S^{-1}J_0F(x)\\
 &=\int_{C_S}g^*(|x/\lambda|)\sqrt{|x/\lambda|}S_\Lambda(\lambda)
[(E_S\frak F_S^t E_S^{-1})J_0F](\lambda)d^\times\lambda\\
 &=\int_{C_S}\sqrt{|{\lambda\over x}|}g(|{\lambda\over x}|)
S_\Lambda(\lambda)[(E_S\frak F_S^t E_S^{-1})J_0F](\lambda)d^\times\lambda\\
&=\int_{C_S}\overline{(E_S\frak F_SE_S^{-1})_\lambda\left(\sqrt{|{\lambda\over x}|}
g(|{\lambda\over x}|)S_\Lambda(\lambda)\right)}(t)J_0F(t)d^\times t
\endaligned $$
where we used the fact that $g$ and $S_\Lambda$ are real-valued functions.

 By [5, Lemma 3.4 and Remark, p. 2467],
 $$\aligned &(e_S\frak F_Se_S^{-1})_\lambda[{1\over\sqrt{|x|}}
g(|{\lambda\over x}|)S_\Lambda(\lambda)](t)
=(\frak F_S)_\lambda[{1\over\sqrt{|x|}}
g(|{\lambda\over x}|)S_\Lambda(\lambda)](t)\\
&=\int_{\Bbb A_S}{1\over\sqrt{|x|}}
g(|{\lambda\over x}|)S_\Lambda(\lambda)\Psi_S(-\lambda t)d\lambda
=\sqrt{|x|}\int_{\Bbb A_S}g(\lambda)S_\Lambda(\lambda x)
\Psi_S(\lambda tx)d\lambda.
\endaligned $$
It follows that
 $$(E_S\frak F_SE_S^{-1})_\lambda\{\sqrt{|{\lambda\over x}|}
g(|{\lambda\over x}|)S_\Lambda(\lambda)\}(t)
=\sqrt{|tx|}\int_{\Bbb A_S}g(\lambda)S_\Lambda(\lambda x)
\Psi_S(\lambda tx)d\lambda.$$
Thus,
$$\aligned V_S(g^*)S_\Lambda E_S \frak F_S^tJE_S^{-1}F(x)
&=\int_{C_S}J_0F(t)d^\times t\sqrt{|tx|}\int_{\Bbb A_S}g(\lambda)
S_\Lambda(\lambda x)\Psi_S(\lambda tx)d\lambda\\
&=\int_{C_S}\left(\sqrt{|x/t|}\int_{\Bbb A_S}g(\lambda)
S_\Lambda(\lambda x)\Psi_S(\lambda x/t)d\lambda\right)F(t)d^\times t.
\endaligned$$
Since $E_S(S(\Bbb A_S)_0)$ is dense in $L^2(C_S)$ and since
$V_S(g^*)S_\Lambda E_S \frak F_S^tJE_S^{-1}$ is bounded on
$L^2(C_S)$, this identity holds for all $F\in L^2(C_S)$.

 This completes the proof of the lemma.
\qed\enddemo

\proclaim{Lemma 3.7}   For each $c\in (0,1)$ we have
$$\left|\int_{\Bbb A_S}g(\lambda)
S_\Lambda(\lambda z)\Psi_S(\lambda {z\over y})d\lambda\right|
\ll_S \left|{z\over y}\right|^{-c}.$$
     \endproclaim

 \demo{Proof}  Since $|\gamma|_S=1$ for $\gamma\in O_S^*$ and since
 the complement of $J_S$ in $\Bbb A_S$ is negligible as $S$ is a finite set,
  by changing $\lambda\to \lambda |z/y|y/z$ also with
  $|z/y|:=(|z/y|,1,\cdots,1)$ depending on context we can write
  $$\aligned &\int_{\Bbb A_S}g(\lambda)
S_\Lambda(\lambda z)\Psi_S(\lambda {z\over y})d\lambda
=\int_{\Bbb A_S}g(\lambda)
S_\Lambda(\lambda z)\Psi_S(\lambda |{z\over y}|)d\lambda\\
&=\sum_{\gamma\in O_S^*}\int_{I_S}g(\lambda)S_\Lambda(\lambda z)
\Psi_S(\gamma\lambda|{z\over y}|)d\lambda\\
&=\sum_{\gamma\in O_S^*}\varpi(\gamma)
\int_0^\infty g(u)S_\Lambda(uz)e^{-2\pi i\gamma u|z/y|}du\\
 &=2\rho_S\sum_{k\in\Bbb N_S}{\mu(k)\over\prod_{p|k}(p-1)}
 \sum_{l\in\Bbb N_S, (l,k)=1}\int_0^\infty
 g(u)S_\Lambda(uz)\cos{2\pi lu|z/y|\over k}du
 \endaligned $$
 where
$\rho_S=\prod_{p\in S^\prime}(1-p^{-1})$ and
 $$\varpi(\gamma)=\prod_{p\in S^\prime} \cases 1-p^{-1}
 &\text{if $|\gamma|_p\leq 1$,}\\
 -p^{-1} &\text{if $|\gamma|_p=p,$} \\
 0 &\text{if $|\gamma|_p>p$}. \endcases\tag 3.2$$

   By [11, Example 10, p. 162],
 $$\int_0^\infty t^{s-1}\cos t\,dt=\Gamma(s)\cos{\pi s\over 2}$$
 for $0<\Re s<1$.  Since $g(u)S_\Lambda(uz)$ as a function of $u$
  has a compact support in $(0,\infty)$,  we have
  $$\aligned &\int_{N2lu\pi |z/y|\over k}^\infty t^{s-1}\cos t\,dt
  =-({N2lu\pi |z/y|\over k})^{s-1}\sin {N2lu\pi |z/y|\over k}-(s-1)
  \int_{N2lu\pi |z/y|\over k}^\infty t^{s-2}\sin t\,dt\\
  &\leq ({N2lu\pi |z/y|\over k})^{\Re s-1}\left(1+{|s-1|\over 1-\Re s}\right)\to 0
  \endaligned $$
 uniformly for $u$ in the support of $g(u)S_\Lambda(uz)$ as $N\to\infty$
 when $0<\Re s<1$.
 This implies that for $0<\Re s<1$
 $$\aligned &\int_0^\infty t^{s-1}dt\int_0^\infty
 g(u)S_\Lambda(uz)\cos{2lu\pi t|z/y|\over k}du\\
 &=\lim_{N\to\infty}\int_{\sqrt\mu}^{1\over\sqrt\mu} g(u)S_\Lambda(uz)du
 \int_0^N t^{s-1}\cos{2lu\pi t|z/y|\over k}dt\\
 &=\lim_{N\to\infty}\int_{\sqrt\mu}^{1\over\sqrt\mu} ({2lu\pi|z/y|\over k})^{-s}
 g(u)S_\Lambda(uz)du
 \int_0^{2lu\pi N|z/y|\over k} t^{s-1}\cos t\,dt\\
 &=k^s(2l\pi|z/y|)^{-s}\int_0^\infty g(u)S_\Lambda(uz)u^{-s}du
 \int_0^\infty t^{s-1}\cos t\,dt\\
 &=k^s(2l\pi|z/y|)^{-s}\widehat{g(u)S_\Lambda(uz)}(1-s)
 \Gamma(s)\cos{\pi s\over 2}
 \endaligned $$
 where the Mellin transform
 $\widehat{g(u)S_\Lambda(uz)}(1-s)=\int_0^\infty g(u)S_\Lambda(uz) u^{-s}du$.

 By Mellin's inversion formula [12, Theorem 28, p. 46],
 $$\int_0^\infty g(u)S_\Lambda(uz)\cos{2\pi lut|{z\over y}|\over k}du
 ={1\over 2\pi i}\int_{c-i\infty}^{c+i\infty} t^{-s}
 k^s(2l\pi|{z\over y}|)^{-s}\widehat{g(u)S_\Lambda(uz)}(1-s)
 \Gamma(s)\cos{\pi s\over 2}ds$$
 for $0<c<1$.  Because $\pm\gamma\in O_S^*$ for every $\gamma\in O_S^*$,
 choosing $t=1$ we get
 $$\aligned &\int_{\Bbb A_S}g(\lambda)
S_\Lambda(\lambda z)\Psi_S(\lambda {z\over y})d\lambda\\
&={\rho_S\over\pi i}\sum_{k\in\Bbb N_S}{\mu(k)\over\prod_{p|k}(p-1)}
 \sum_{l\in\Bbb N_S, (l,k)=1}\int_{c-i\infty}^{c+i\infty}
 k^s(2l\pi|{z\over y}|)^{-s}\widehat{g(u)S_\Lambda(uz)}(1-s)
 \Gamma(s)\cos{\pi s\over 2}ds\\
 &={1\over 2\pi i}\int_{c-i\infty}^{c+i\infty}
    |z/y|^{-s}\widehat{g(u)S_\Lambda(uz)}(1-s)\chi(1-s)
    \prod_{p\in S^\prime}{1-p^{s-1}\over 1-p^{-s}}ds\endaligned$$
 where $\chi(s)=2^s\pi^{s-1}\Gamma(1-s)\sin{\pi s\over 2}$.
 By partial integration,
 $$\aligned &\widehat{g(u)S_\Lambda(uz)}(1-s)
 =\int_{\max(1,{1\over\Lambda|z|})}^{1\over\mu}g(u)u^{-s}du
 ={1\over s-1}\{g(\max(1,{1\over\Lambda|z|}))
 (\max(1,{1\over\Lambda|z|}))^{1-s}\\
 &-g({1\over\mu})\mu^{s-1}
 -\int_0^\infty g^\prime(u)S_\Lambda(uz)u^{1-s}du\}\ll {1\over |s-1|}.
 \endaligned $$
  Stirling's formula (cf. [11, line 13, p.151]) reads that
  for any fixed value of $\sigma$
  $$|\Gamma(\sigma+it)|\sim \sqrt{2\pi}e^{-{\pi\over 2}|t|}|t|^{\sigma-{1\over 2}}$$
  as $|t|\to\infty$.  Thus, we can move line of integration so that
  $$\left|\int_{\Bbb A_S}g(\lambda)
S_\Lambda(\lambda z)\Psi_S(\lambda {z\over y})d\lambda\right|\ll_S |z/y|^{-c}$$
for each $0<c<1$.

 This completes the proof of the lemma.
  \qed\enddemo

\demo{Proof of Theorem 1.3}  By Lemma 3.2,
$$V(h)T=-V(g)[S_\Lambda, V(g^*)]
 S_\Lambda+V(g)[S_\Lambda,
 V(g^*)S_\Lambda E_S \frak F_S^tJE_S^{-1}]
E_S J\frak F_S E_S^{-1} S_\Lambda.\tag 3.3$$
Since $E_S JE_S^{-1}F(x)=F(1/x)$ by (3.1) and since
$E_S J\frak F_S E_S^{-1}=E_S JE_S^{-1}\cdot E_S\frak F_S E_S^{-1}$,
both $S_\Lambda$ and $E_S J\frak F_S E_S^{-1}$ are bounded operators
on $L^2(C_S)$.  By Lemma 3.4, in order to show that $V(h)T$ is of trace class
on $L^2(C_S)$ it suffices to show that
$V(g)[S_\Lambda, V(g^*)]$ and
$$V(g)[S_\Lambda, V(g^*)S_\Lambda E_S \frak F_S^tJE_S^{-1}]$$
 are of trace class.

We can write
$$\aligned  &V_S(g)[S_\Lambda, V_S(g^*)]F(x)\\
&=\int_{C_S}\sqrt{|x/z|}g(|x/z|)d^\times z\int_{C_S}
 \left(S_\Lambda(|z|)-S_\Lambda(|y|)\right) g^*(|z/y|)
 \sqrt{|z/y|}F(y)d^\times y\\
 &=\int_{C_S}{g(|{x\over z}|)\sqrt{|{x\over z}|}\over 1+|\log|z||}d^\times z
 \int_{C_S}  {[S_\Lambda(|z|)-S_\Lambda(|y|)]
 g^*(|{z\over y}|)\sqrt{|{z\over y}|}\over (1+|\log|z||)^{-1}}F(y)d^\times y
 =ABF(x)
\endaligned \tag 3.4$$
for $F\in L^2(C_S)$, where the operators $A$ and $B$ are defined by
$$BF(z)=\int_{C_S}  {[S_\Lambda(|z|)-S_\Lambda(|y|)]
 g^*(|{z\over y}|)\sqrt{|{z\over y}|}\over (1+|\log|z||)^{-1}}F(y)d^\times y$$
 and
 $$AG(x)=\int_{C_S}{g(|{x\over z}|)\sqrt{|{x\over z}|}\over 1+|\log|z||}
 G(z)d^\times z.$$

  We find that
$$\aligned
&\int_{C_S\times C_S}\left|{g
\left(|{x\over z}|\right)\sqrt{|{x\over z}|}\over 1+|\log|z||}\right|^2d^\times xd^\times z\\
&=\int_{C_S}{d^\times z\over (1+|\log|z||)^2}
\int_{C_S}\left|g\left(\left|{x\over z}\right|\right)
\left|{x\over z}\right|^{1\over 2}\right|^2d^\times x\\
&=\int_0^\infty {d^\times |z|\over (1+|\log|z||)^2}\int_0^\infty
|g(|x|)\sqrt{|x|}|^2d^\times |x|<\infty
 \endaligned $$
as $g(x)$ has a compact support.  By Lemma 3.5, $A$ is a Hilbert-Schmidt operator.

 If
 $$[S_\Lambda(|y|)-S_\Lambda(|z|)]g^*(|{z\over y}|)\neq 0,$$
 then we must have $\sqrt\mu\leq |z/y|\leq \mu^{-1/2}$ and
 either $|y|>\Lambda^{-1}$, $|z|\leq\Lambda^{-1}$
 or $|z|>\Lambda^{-1}$, $|y|\leq\Lambda^{-1}$.

 If $|y|>\Lambda^{-1}$ and $|z|\leq\Lambda^{-1}$, then
 $$\Lambda^{-1}<|y|\leq |z|/\sqrt\mu\leq {1\over\sqrt\mu\Lambda}$$
  and
$$\Lambda^{-1}\geq |z|\geq |y|\sqrt\mu>{\sqrt\mu\over\Lambda}.$$

If $|z|>\Lambda^{-1}$, $|y|\leq\Lambda^{-1}$, then
$$\Lambda^{-1}<|z|\leq |y|/\sqrt\mu\leq {1\over\sqrt\mu\Lambda}$$
and
$$\Lambda^{-1}\geq |y|\geq \sqrt\mu |z|>{\sqrt\mu\over\Lambda}.$$
Therefore,
 $[S_\Lambda(|y|)-S_\Lambda(|z|)]g^*(|{z\over y}|)\neq 0$ only if either
 $$|z|\in ({\sqrt\mu\over\Lambda}, {1\over\Lambda}],\,\,
|y|\in ({1\over \Lambda}, {1\over\Lambda\sqrt\mu}]\,\,\,\, \text{ or }\,\,\,\,
|z|\in ({1\over \Lambda}, {1\over\Lambda\sqrt\mu}],\,\,
|y|\in ({\sqrt\mu\over\Lambda}, {1\over\Lambda}].\tag 3.5 $$
Hence
$$\aligned
&\int_{C_S\times C_S}\left|(1+|\log|z||)[S_\Lambda(|z|)-S_\Lambda(|y|)]
 g^*(|{z\over y}|)\sqrt{|{z\over y}|}\right|^2d^\times yd^\times z\\
 &\leq \int_{\sqrt\mu\over\Lambda}^{1\over\Lambda\sqrt\mu}
 (1+|\log|z||)^2d^\times z\int_{C_S}
 |g^*(|z/y|)|^2|z/y|d^\times y\\
 &=\int_{\sqrt\mu\over\Lambda}^{1\over\Lambda\sqrt\mu}
 (1+|\log|z||)^2d^\times z\int_{C_S}
 |g^*(|y|)|^2|y|d^\times y<\infty.
\endaligned $$
It follows from Lemma 3.5 that $B$ is a Hilbert-Schmidt operator.

Thus $V(g)[S_\Lambda, V(g^*)]$ is a product of two
Hilbert-Schmidt operators $A$, $B$ on $L^2(C_S)$, and hence is of trace class
on $L^2(C_S)$ by Lemma 3.3.

By Lemma 3.6, we can write
$$\aligned  &V(g)[S_\Lambda, V(g^*)S_\Lambda E_S
\frak F_S^tJE_S^{-1}]F(x)=\int_{C_S}{g(|{x\over z}|)
\sqrt{|{x\over z}|}\over 1+|\log|z||}d^\times z\\
&\times \int_{C_S} {S_\Lambda(|z|)-S_\Lambda(|y|)\over (1+|\log|z||)^{-1}}
\left(\sqrt{|{z\over y}|}\int_{\Bbb A_S}g(\lambda)
S_\Lambda(\lambda z)\Psi_S(\lambda {z\over y})d\lambda\right)F(y)d^\times y
\endaligned \tag 3.6$$
for $F\in L^2(C_S)$.  We have already shown above that
$$\int_{C_S\times C_S}\left|{g\left(|{x\over z}|\right)
\sqrt{|{x\over z}|}\over 1+|\log|z||}\right|^2d^\times xd^\times z<\infty. $$

Since $S_\Lambda(|z|)-S_\Lambda(|y|)\neq 0$ only if either $\Lambda |z|>1$,
$\Lambda|y|\leq 1$ or $\Lambda|y|>1$, $\Lambda |z|\leq 1$.
That is,
$$|\log|z/y||=|\log|\Lambda z/\Lambda y||
=|\log|\Lambda z|-\log|\Lambda y||
=|\log|\Lambda z||+|\log|\Lambda y||.\tag 3.7$$
By Lemma 3.7
$$\sqrt{|{z\over y}|}\int_{\Bbb A_S}g(\lambda)
S_\Lambda(\lambda z)\Psi_S(\lambda {z\over y})d\lambda
\ll_S \left|{\Lambda z\over \Lambda y}\right|^{{1\over 2}-c}.$$
If we choose $c=1/4$ when $|\Lambda z/\Lambda y|<1$ and choose $c=3/4$
when $|\Lambda z/\Lambda y|>1$, by (3.7)
$$\sqrt{|{z\over y}|}\left|\int_{\Bbb A_S}g(\lambda)
S_\Lambda(\lambda z)\Psi_S(\lambda {z\over y})d\lambda\right|\ll_S
{1\over (2+|\log|z/y||)^3}={1\over (1+|\log|\Lambda z||+1+|\log|\Lambda y||)^3}$$
for all $y,z\in C_S$ satisfying $S_\Lambda(|z|)-S_\Lambda(|y|)\neq 0$.
It follows that
$$\aligned
&\int_{C_S\times C_S}\left|{S_\Lambda(|z|)-S_\Lambda(|y|)\over (1+|\log|z||)^{-1}}
\left(\sqrt{|{z\over y}|}\int_{\Bbb A_S}g(\lambda)
S_\Lambda(\lambda z)\Psi_S(\lambda {z\over y})d\lambda\right)\right|^2
d^\times yd^\times z\\
 &\ll_S \int_{C_S}  {(1+|\log|z||)^2\over (1+|\log|\Lambda z||)^4}d^\times z
 \int_{C_S} {1\over (1+|\log|\Lambda y||)^2}d^\times y<\infty.
 \endaligned \tag 3.8$$
Thus $V(g)[S_\Lambda, V(g^*)S_\Lambda E_S\frak F_S^tJE_S^{-1}]$
is a product of two Hilbert-Schmidt operators by Lemma 3.5, and
hence is of trace class on $L^2(C_S)$ by Lemma 3.3.

 This completes the proof of Theorem 1.3.
\qed\enddemo

\vskip 0.35in
\heading
4.  Trace of $V(h)T$
\endheading

  In this section we compute the trace of $V(h)T$ on $L^2(C_S)$.
    The $S$-local trace formula is first given in
  [3, Theorem 4, p. 56] for the operator
  $\frak F_SP_\Lambda\frak F_S^t P_\Lambda$.
  The author's variant of it in  [5, Theorem 3.16, p. 2481]
  is due to Meyer [7, (19), p. 549]
  (although the author discovered it independently without
  knowing its existence already in [7]).

\proclaim{Lemma 4.1} ([2, Corollary 3.2, p. 237])  Let $\mu$ be a $\sigma$-finite
Borel measure on a second countable space $X$, and  let
 $A$ be a trace class Hilbert-Schmidt integral operator on
$L^2(X, d\mu)$.  If the kernel $k(x, y)$ is continuous at $(x,x)$
for almost every $x$, then
$$\text{trace}(A)=\int_Xk(x,x)d\mu(x).$$
\endproclaim

\proclaim{Lemma 4.2}  $V(h)T$ is a Hilbert-Schmidt integral operator on $L^2(C_S)$ with
 the kernel
 $$\aligned &k(x,y)=-\int_{C_S}{g(|{x\over z}|)
\sqrt{|{x\over z}|}\over 1+|\log|z||}\cdot
 {[S_\Lambda(|z|)-S_\Lambda(|y|)] g^*(|{z\over y}|)
 \sqrt{|{z\over y}|}\over (1+|\log|z||)^{-1}}S_\Lambda(y)d^\times z
 +\int_{C_S}g(|{x\over z}|)\sqrt{|{x\over z}|}\\
 &\times\left(E_S\frak F_S E_S^{-1}\right)_u\left[[S_\Lambda(z)
-S_\Lambda({1\over |u|})]\sqrt{|zu|}\int_{\Bbb A_S}g(\lambda)
S_\Lambda(\lambda z)\Psi_S(\lambda zu)d\lambda\right](y)
S_\Lambda(y)d^\times z.
 \endaligned $$
\endproclaim

\demo{Proof}  By (3.3) and (3.4),
$$V(h)T=-V(g)[S_\Lambda, V(g^*)]
 S_\Lambda+V(g)[S_\Lambda,
 V(g^*)S_\Lambda E_S \frak F_S^tJE_S^{-1}]
E_S J\frak F_S E_S^{-1} S_\Lambda \tag 4.1$$
and
$$\aligned &V(g)[S_\Lambda, V(g^*)]S_\Lambda F(x)\\
 &=\int_{C_S}{g(|{x\over z}|)\sqrt{|{x\over z}|}\over 1+|\log|z||}d^\times z
 \int_{C_S}  {[S_\Lambda(|z|)-S_\Lambda(|y|)] g^*(|{z\over y}|)
 \sqrt{|{z\over y}|}\over (1+|\log|z||)^{-1}}S_\Lambda(y)F(y)d^\times y.
 \endaligned\tag 4.2$$
 Since
 $$E_SJf(y)=\sqrt{|y|}\sum_{\xi\in O_S^*}{1\over|\xi y|}f({1\over\xi y})
 ={1\over\sqrt{|y|}}\sum_{\gamma\in O_S^*}f({\gamma\over y})=E_S(f)({1\over y}),$$
  by (3.6) we have
 $$\aligned  &V(g)[S_\Lambda, V(g^*)S_\Lambda E_S
\frak F_S^tJE_S^{-1}]E_S J\frak F_S E_S^{-1} S_\Lambda F(x)
=\int_{C_S}g(|{x\over z}|)\sqrt{|{x\over z}|}d^\times z\\
&\times\int_{C_S} [S_\Lambda(|z|)-S_\Lambda(|y|)]
\{\sqrt{|{z\over y}|}\int_{\Bbb A_S}g(\lambda)
S_\Lambda(\lambda z)\Psi_S(\lambda {z\over y})d\lambda\}
\left(E_S J\frak F_S E_S^{-1} S_\Lambda F\right)(y)d^\times y\\
&=\int_{C_S}g(|{x\over z}|)\sqrt{|{x\over z}|}d^\times z
 \int_{C_S} [S_\Lambda(|z|)-S_\Lambda(|y|)]\\
 &\times\{\sqrt{|{z\over y}|}\int_{\Bbb A_S}g(\lambda)
S_\Lambda(\lambda z)\Psi_S(\lambda {z\over y})d\lambda\}
\left(E_S\frak F_S E_S^{-1} S_\Lambda F\right)({1\over y})d^\times y\\
&=\int_{C_S}g(|{x\over z}|)\sqrt{|{x\over z}|}d^\times z
\int_{C_S} [S_\Lambda(|z|)-S_\Lambda({1\over |y|})]\\
&\times\{\sqrt{|zy|}\int_{\Bbb A_S}g(\lambda)
S_\Lambda(\lambda z)\Psi_S(\lambda zy)d\lambda\}
\left(E_S\frak F_S E_S^{-1} S_\Lambda F\right)(y)d^\times y
\endaligned $$
 for $F\in L^2(C_S)$.

 Since
$$[S_\Lambda(|z|)-S_\Lambda(|y|)]
\{\sqrt{|{z\over y}|}\int_{\Bbb A_S}g(\lambda)
S_\Lambda(\lambda z)\Psi_S(\lambda {z\over y})d\lambda\}\in L^2(C_S)$$
for fixed $z$ by (3.8), we have
$$[S_\Lambda(|z|)-S_\Lambda({1\over |y|})]
\{\sqrt{|zy|}\int_{\Bbb A_S}g(\lambda)
S_\Lambda(\lambda z)\Psi_S(\lambda zy)d\lambda\}\in L^2(C_S)$$
as a function of $y$. By Lemma 2.1, $E_S\frak F_S E_S^{-1}$
is a unitary operator on $L^2(C_S)$.  This implies that
$$\langle E_S\frak F_S E_S^{-1}F, G\rangle_{L^2(c_S)}
=\langle F, E_S\frak F_S^t E_S^{-1}G\rangle_{L^2(C_S)}$$
for $F, G\in L^2(C_S)$.  From this identity we derive that
$$\aligned &\int_{C_S} [S_\Lambda(|z|)
-S_\Lambda({1\over |y|})]\{\sqrt{|zy|}\int_{\Bbb A_S}g(\lambda)
S_\Lambda(\lambda z)\Psi_S(\lambda zy)d\lambda\}
\left(E_S\frak F_S^t E_S^{-1} S_\Lambda F\right)(y)d^\times y\\
&=\int_{C_S} \left(E_S\frak F_S E_S^{-1}\right)_u\left[[S_\Lambda(z)
-S_\Lambda({1\over |u|})]\sqrt{|zu|}\int_{\Bbb A_S}g(\lambda)
S_\Lambda(\lambda z)\Psi_S(\lambda zu)d\lambda\right](y)
S_\Lambda F(y)d^\times y.
\endaligned $$
Hence,
 $$\aligned  &V(g)[S_\Lambda, V(g^*)S_\Lambda E_S
\frak F_S^tJE_S^{-1}]E_S J\frak F_S E_S^{-1} S_\Lambda F(x)
=\int_{C_S}g(|{x\over z}|)\sqrt{|{x\over z}|}d^\times z\\
&\times\int_{C_S} \left(E_S\frak F_S E_S^{-1}\right)_u\left[[S_\Lambda(z)
-S_\Lambda({1\over |u|})]\sqrt{|zu|}\int_{\Bbb A_S}g(\lambda)
S_\Lambda(\lambda z)\Psi_S(\lambda zu)d\lambda\right](y)
S_\Lambda F(y)d^\times y
\endaligned\tag 4.3 $$
 for $F\in L^2(C_S)$.

Let
$$\aligned &k(x,y)=-\int_{C_S}{g(|{x\over z}|)
\sqrt{|{x\over z}|}\over 1+|\log|z||}\cdot
 {[S_\Lambda(|z|)-S_\Lambda(|y|)] g^*(|{z\over y}|)
 \sqrt{|{z\over y}|}\over (1+|\log|z||)^{-1}}S_\Lambda(y)d^\times z
 +\int_{C_S}g(|{x\over z}|)\sqrt{|{x\over z}|}\\
 &\times\left(E_S\frak F_S E_S^{-1}\right)_u\left[[S_\Lambda(z)
-S_\Lambda({1\over |u|})]\sqrt{|zu|}\int_{\Bbb A_S}g(\lambda)
S_\Lambda(\lambda z)\Psi_S(\lambda zu)d\lambda\right](y)
S_\Lambda(y)d^\times z.
 \endaligned$$
 Then
 $$ V(h)TF(x)=\int_{C_S}k(x,y)F(y)d^\times y,$$
for $F\in L^2(C_S)$, where the change order of integration
is permissible as the involved integrals are absolutely integrable by
(4.1), (4.2), (4.3).  More precisely,  by (3.5) the double integral
in (4.2) is absolute integrable as $g(t)=0$ when
$t\not\in (\sqrt\mu, \mu^{-1/2})$.  By Schwarz inequality, Lemma 2.1, and (3.8),
we have
$$\aligned &\int_{C_S}| \left(E_S\frak F_S E_S^{-1}\right)_u\left[[S_\Lambda(z)
-S_\Lambda({1\over |u|})]\sqrt{|zu|}\int_{\Bbb A_S}g(\lambda)
S_\Lambda(\lambda z)\Psi_S(\lambda zu)d\lambda\right](y)
S_\Lambda F(y)|d^\times y\\
&\leq (\int_{C_S}|[S_\Lambda(z)-S_\Lambda(y)]
\sqrt{|{z\over y}|}\int_{\Bbb A_S}g(\lambda)S_\Lambda(\lambda z)
\Psi_S(\lambda {z\over y})d\lambda|^2d^\times y)^{1\over 2}
(\int_{C_S}|S_\Lambda F(y)|^2d^\times y)^{1\over 2}<\infty.
\endaligned $$
This implies that the double integral in (4.3) is absolutely integrable.

Next we show that $V(h)T$ is a Hilbert-Schmidt integral operator
on $L^2(C_S)$, i.e., $k(x,y)\in L^2(C_S\times C_S)$.  Notice that
 $$\aligned &\int_{C_S\times C_S}|k(x,y)|^2d^\times xd^\times y
\leq 2 \int_{C_S\times C_S}\left|{g(|{x\over z}|)
\sqrt{|{x\over z}|}\over 1+|\log|z||}\right|^2d^\times zd^\times x\\
&\times\{\int_{C_S\times C_S}\left| {[S_\Lambda(|z|)-S_\Lambda(|y|)]
g^*(|{z\over y}|)\sqrt{|{z\over y}|}
\over (1+|\log|z||)^{-1}}S_\Lambda(y)\right|^2 d^\times zd^\times y
+\int_{C_S\times C_S}\\
&\left|{\left(E_S\frak F_S E_S^{-1}\right)_u
\left[[S_\Lambda(z)-S_\Lambda({1\over |u|})]\sqrt{|zu|}
\int_{\Bbb A_S}g(\lambda)S_\Lambda(\lambda z)\Psi_S(\lambda zu)
d\lambda\right](y)S_\Lambda(y)\over(1+|\log|z||)^{-1}}\right|^2
 d^\times zd^\times y\}.
\endaligned$$
By Lemma 2.1, $E_S\frak F_S^tE_S^{-1}$ is unitary on $L^2(C_S)$.  Hence,
$$\aligned &\int_{C_S}\left|\left(E_S\frak F_S E_S^{-1}\right)_u
\left[[S_\Lambda(z)-S_\Lambda({1\over |u|})]\sqrt{|zu|}
\int_{\Bbb A_S}g(\lambda)S_\Lambda(\lambda z)\Psi_S(\lambda zu)
d\lambda\right](y)S_\Lambda(y)\right|^2d^\times y\\
&\leq \int_{C_S}\left|\left(E_S\frak F_S E_S^{-1}\right)_u
\left[[S_\Lambda(z)-S_\Lambda({1\over |u|})]\sqrt{|zu|}
\int_{\Bbb A_S}g(\lambda)S_\Lambda(\lambda z)\Psi_S(\lambda zu)
d\lambda\right](y)\right|^2d^\times y\\
 &=\int_{C_S}\left|[S_\Lambda(z)-S_\Lambda({1\over |y|})]\sqrt{|zy|}
\int_{\Bbb A_S}g(\lambda)S_\Lambda(\lambda z)\Psi_S(\lambda zy)
d\lambda\right|^2d^\times y\\
 &=\int_{C_S}\left|[S_\Lambda(z)-S_\Lambda(|y|)]\sqrt{|z/y|}
\int_{\Bbb A_S}g(\lambda)S_\Lambda(\lambda z)\Psi_S(\lambda z/y)
d\lambda\right|^2d^\times y.
 \endaligned $$
 By (3.5) and (3.8),
 $$\aligned &\int_{C_S\times C_S}|k(x,y)|^2d^\times xd^\times y
\leq 2 \int_{C_S\times C_S}\left|{g(|{x\over z}|)
\sqrt{|{x\over z}|}\over 1+|\log|z||}\right|^2d^\times zd^\times x\\
&\times\{\int_{C_S\times C_S}\left| {[S_\Lambda(z)-S_\Lambda(y)]
g^*(|{z\over y}|)\sqrt{|{z\over y}|}
\over (1+|\log|z||)^{-1}}S_\Lambda(y)\right|^2 d^\times zd^\times y
+\int_{C_S\times C_S}\\
&\left|{[S_\Lambda(z)-S_\Lambda(y)]\sqrt{|z/y|}
\int_{\Bbb A_S}g(\lambda)S_\Lambda(\lambda z)\Psi_S(\lambda z/y)
d\lambda\over(1+|\log|z||)^{-1}}\right|^2
 d^\times zd^\times y\}<\infty.
\endaligned$$
Therefore $k(x,y)\in L^2(C_S\times C_S)$, and hence
$V(h)T$ is a Hilbert-Schmidt operator on $L^2(C_S)$.

  This completes the proof of the lemma.
\qed\enddemo

\proclaim{Lemma 4.3}  Let $O_p=\{x\in Q_p: |x|_p\leq 1\}$ and
$B=\prod_{p\in S^\prime}p^{-1}O_p$. If $F(y)=f(|y|)$ for some
function $f$, then
$$\int_{\Bbb A_S}F(y)\Psi_S(yx)dy=\int_{\Bbb R\times B}F(y)\Psi_S(y|x|)dy.$$
\endproclaim

\demo{Proof} Since the complement of $J_S$ in $\Bbb A_S$ is negligible,
as $|\xi|=1$ for $\xi\in O_S^*$ we can write
$$\int_{\Bbb A_S}F(y)\Psi_S(yx)dy
=\sum_{\xi\in O_S^*}\int_{I_S}F(y)\Psi_S(yx\xi)dy.$$
Since $F(y)=f(|y|)$, if we change variables $y\to y|x|/x$
with $|x|=(|x|, 1, \cdots,1)$ then
$$\int_{\Bbb A_S}F(y)\Psi_S(yx)dy
=\sum_{\xi\in O_S^*}\int_{I_S}f(|y|)\Psi_S(y|x|\xi)dy
=\sum_{\xi\in O_S^*}\varpi(\xi)\int_0^\infty f(t)e^{-2\pi it|x|\xi}dt,$$
where $\varpi(\xi)$ is given as in (3.2).

  If $\xi\in O_S^*$ satisfies $\varpi(\xi)\neq 0$,
   then $\xi=l/k$ for some integers $k,l\in\Bbb N_S$
   with square free $k$.  Thus, $\xi\prod_{p\in S^\prime}O_p^*\subset B$. Hence,
   $$\bigcup_{\xi\in O_S^*, \varpi(\xi)\neq 0}\xi I_S\subset\Bbb R^*\times B.$$
   Conversely, if $x=(x_p)\in\Bbb R^*\times B$ then there exists a $\xi\in O_S^*$
   such that $\xi^{-1}x\in I_S$, i.e., $|\xi|_p=|x|_p\leq p$ for $p\in S^\prime$.
      By (3.2), $\varpi(\xi)\neq 0$.  Hence,
   $x\in \bigcup_{\xi\in O_S^*, \varpi(\xi)\neq 0}\xi I_S$.  Thus,
   $$\Bbb R^*\times B\subset \bigcup_{\xi\in O_S^*, \varpi(\xi)\neq 0}\xi I_S.$$
     Also if $\xi, \gamma$ are distinct elements in $O_S^*$, then
  $\xi O_S^*$ and $\gamma O_S^*$ are disjoint sets.  Therefore,
  $$\Bbb R^*\times B=\bigcup_{\xi\in O_S^*, \varpi(\xi)\neq 0}\xi I_S,
  \text{ a disjoint union}.$$
  Therefore,
  $$\sum_{\xi\in O_S^*, \varpi(\xi)\neq 0}\int_{I_S}F(y)\Psi_S(y|x|\xi)dy
  =\int_{\Bbb R\times B}F(y)\Psi_S(y|x|\xi)dy.$$
  That is,
  $$\int_{\Bbb A_S}F(y)\Psi_S(yx)dy=\int_{\Bbb R\times B}F(y)\Psi_S(y|x|)dy.$$

This completes the proof of the lemma.
\qed\enddemo

\proclaim{Lemma 4.4} We have
$$\left(E_S\frak F_S E_S^{-1}\right)_u
\left[\sqrt{|u|}\int_{\Bbb A_S}g(\lambda)
S_\Lambda(\lambda z)\Psi_S(\lambda zu)d\lambda\right](y)
 ={\sqrt{|y|}\over |z|}g(|{y\over z}|)S_\Lambda(y).$$
That is,
 $$\int_{\Bbb A_S}\Psi_S(-uy)du
\int_{\Bbb A_S}{1\over |z|}g({\lambda\over z})
S_\Lambda(\lambda)\Psi_S(\lambda u)d\lambda
={1\over |z|}g({y\over z})S_\Lambda(y).$$
\endproclaim

\demo{Proof} By [5, Lemma 3.4 and Remark, p. 2467],
 $$\left(e_S\frak F_S e_S^{-1}\right)_u[\int_{\Bbb A_S}g(\lambda)
S_\Lambda(\lambda z)\Psi_S(\lambda zu)d\lambda](y)
=\left(\frak F_S\right)_u[\int_{\Bbb A_S}g(\lambda)
S_\Lambda(\lambda z)\Psi_S(\lambda zu)d\lambda](y).$$
It follows that
  $$\aligned &\left(E_S\frak F_S E_S^{-1}\right)_u
  [\sqrt{|u|}\int_{\Bbb A_S}g(\lambda)
S_\Lambda(\lambda z)\Psi_S(\lambda zu)d\lambda](y)\\
&=\sqrt{|y|}\left(\frak F_S\right)_u[\int_{\Bbb A_S}g(\lambda)
S_\Lambda(\lambda z)\Psi_S(\lambda zu)d\lambda](y)
=\sqrt{|y|}\left(\frak F_S\right)_u[\int_{\Bbb A_S}{1\over |z|}g({\lambda\over z})
S_\Lambda(\lambda)\Psi_S(\lambda u)d\lambda](y)\\
&={\sqrt{|y|}\over|z|}\int_{\Bbb A_S}\Psi_S(-uy)du\int_{\Bbb A_S}
g({\lambda\over z})S_\Lambda(\lambda)\Psi_S(\lambda u)d\lambda.
\endaligned$$
Hence, the second stated identity follows from the first one.

 By Lemma 4.3,
$$\aligned &\int_{\Bbb A_S}\Psi_S(-uy)du\int_{\Bbb A_S}
g({\lambda\over z})S_\Lambda(\lambda)\Psi_S(\lambda u)d\lambda\\
&=\int_{\Bbb A_S}\Psi_S(-uy)du\int_{\Bbb R\times B}
g({\lambda\over z})S_\Lambda(\lambda)\Psi_S(\lambda|u|)d\lambda\\
&=\int_{\Bbb R\times B}\Psi_S(-uy)du\int_{\Bbb R\times B}
g({\lambda\over z})S_\Lambda(\lambda)\Psi_S(\lambda u)d\lambda.
\endaligned$$
 Denote
  $$x=(x_r, x_b)\text{ with } x_r\in\Bbb R, x_b\in B. \tag 4.4$$
We can write
$$\aligned &\int_{\Bbb A_S}\Psi_S(-uy)du\int_{\Bbb A_S}
g({\lambda\over z})S_\Lambda(\lambda)\Psi_S(\lambda u)d\lambda\\
&=\int_{-\infty}^\infty e^{2\pi iu_ry_r}du_r \int_{\Bbb A_{S^\prime}}
\Psi_{S^\prime}(-u_by_b)du_b\int_{\Bbb A_{S^\prime}}\phi(\lambda_b)
\Psi_{S^\prime}(\lambda_bu_b)d\lambda_b \endaligned $$
where
$$\phi(\lambda_b)=\cases \int_{-\infty}^\infty
g({|\lambda_r||\lambda_b|\over|z|})S_\Lambda(|\lambda_r||\lambda_b|)
 e^{-2\pi i\lambda_ru_r}d\lambda_r &\text{if $\lambda_b\in B$}\\
0&\text{if $\lambda_b\not\in B$.}\endcases$$
Since $\phi(\lambda_b)$ is locally constant as
a function of $\lambda_b\in\Bbb A_{S^\prime}$ and is supported on the compact set
$B$, the condition of [10, Theorem 2.2.2, p. 310] is satisfied by this function.
By the Fourier inversion formula (Note that [10, Theorem 2.2.2, p. 310] is still
true if we replace $k_p^+$ there by $\Bbb A_{S^\prime}$),  we have
$$\aligned &\int_{\Bbb A_{S^\prime}}
\Psi_{S^\prime}(-u_by_b)du_b\int_{\Bbb A_{S^\prime}}\phi(\lambda_b)
\Psi_{S^\prime}(\lambda_bu_b)d\lambda_b=\phi(y_b)\\
&=\int_{-\infty}^\infty
g({|\lambda_r||y_b|\over|z|})S_\Lambda(|\lambda_r||y_b|)
 e^{-2\pi i\lambda_ru_r}d\lambda_r.\endaligned$$
 Since $g({t|y_b|\over|z|})S_\Lambda(t|y_b|)$
 is a continuous and compactly supported
function of $t\in(0,\infty)$, $t\neq 1/\Lambda|y_b|$
 and is of bounded variation in an interval including
$y_r$, by Fourier's single-integral formula [12, Theorem 12, p. 25] we have
$$\aligned &\int_{\Bbb A_S}\Psi_S(-uy)du\int_{\Bbb A_S}
g({\lambda\over z})S_\Lambda(\lambda)\Psi_S(\lambda u)d\lambda\\
&=\int_{-\infty}^\infty e^{2\pi iu_ry_r}du_r \int_{-\infty}^\infty
g({|\lambda_r||y_b|\over|z|})S_\Lambda(|\lambda_r||y_b|)
 e^{-2\pi i\lambda_ru_r}d\lambda_r\\
 &=\int_{-\infty}^\infty e^{2\pi iu_ry_r}du_r \int_{|z|\sqrt\mu\over|y_b|}
 ^{|z|\over\sqrt\mu|y_b|}g({|\lambda_r||y_b|\over|z|})S_\Lambda(|\lambda_r||y_b|)
 e^{-2\pi i\lambda_ru_r}d\lambda_r\\
 &=\lim_{X\to\infty}\int_{|z|\sqrt\mu\over|y_b|}
 ^{|z|\over\sqrt\mu|y_b|}g({|\lambda_r||y_b|\over|z|})S_\Lambda(|\lambda_r||y_b|)
 d\lambda_r\int_{-X}^X e^{2\pi iu_r(y_r-\lambda_r)}du_r\\
 &=\lim_{X\to\infty}{1\over\pi}\int_{|z|\sqrt\mu\over|y_b|}
 ^{|z|\over\sqrt\mu|y_b|}g({|\lambda_r||y_b|\over|z|})S_\Lambda(|\lambda_r||y_b|)
  {\sin 2\pi X(y_r-\lambda_r)\over y_r-\lambda_r}d\lambda_r
  =g(|y/z|)S_\Lambda(y).
\endaligned$$
Therefore,
$$\left(E_S\frak F_S E_S^{-1}\right)_u
  [\sqrt{|u|}\int_{\Bbb A_S}g(\lambda)
S_\Lambda(\lambda z)\Psi_S(\lambda zu)d\lambda](y)
={\sqrt{|y|}\over|z|}g(|{y\over z}|)S_\Lambda(y).$$

This completes the proof of the lemma.
\qed\enddemo

\proclaim{Lemma 4.5}  Let $k(x,y)$ be given as in Lemma 4.2.  Then
$$k(x,x)=S_\Lambda(x)\int_{\Bbb A_S, |u|\geq\Lambda |x|}\Psi_S(-u)du
\int_{\Bbb A_S}h(\lambda)S_\Lambda(\lambda x)\Psi_S(\lambda u)d\lambda.$$
\endproclaim

\demo{Proof}  By Lemma 4.2,
$$\aligned &k(x,y)=S_\Lambda(y)\sqrt{|x|}
\{-{1\over\sqrt{|y|}}\int_{C_S}g(|{x\over z}|)
[S_\Lambda(|z|)-S_\Lambda(|y|)] g^*(|{z\over y}|)d^\times z\\
 &+\int_{C_S}g(|{x\over z}|)
 \left(E_S\frak F_S E_S^{-1}\right)_u\left[[S_\Lambda(z)
-S_\Lambda({1\over |u|})]\sqrt{|u|}\int_{\Bbb A_S}g(\lambda)
S_\Lambda(\lambda z)\Psi_S(\lambda zu)d\lambda\right](y)d^\times z\}.
 \endaligned $$

Notice that
$$\aligned \int_{C_S}g(|{x\over z}|)g^*(|{z\over y}|)d^\times z
&=\int_0^\infty g(|xz|)|y|g(|yz|)d|z|\\
&=\int_0^\infty g(|{x\over y}z|)g(|z|)d|z|=h(x/y).
\endaligned$$
 By Lemma 4.4,
$$\left(E_S\frak F_S E_S^{-1}\right)_u
\left[\sqrt{|u|}\int_{\Bbb A_S}g(\lambda)
S_\Lambda(\lambda z)\Psi_S(\lambda zu)d\lambda\right](y)
 ={\sqrt{|y|}\over |z|}g(|y/z|)S_\Lambda(y)$$
 for $y\in C_S$.  Hence,
 $$\aligned &\int_{C_S}g(|{x\over z}|)
 \left(E_S\frak F_S E_S^{-1}\right)_u\left[S_\Lambda(z)
 \sqrt{|u|}\int_{\Bbb A_S}g(\lambda)
S_\Lambda(\lambda z)\Psi_S(\lambda zu)d\lambda\right](y)d^\times z\\
&=\sqrt{|y|}S_\Lambda(y)\int_{C_S}g(|{x\over z}|)
S_\Lambda(z){1\over|z|}g(|{y\over z}|)d^\times z\\
&={S_\Lambda(y)\over\sqrt{|y|}}\int_{C_S}g(|{x\over z}|)S_\Lambda(z)
g^*(|{z\over y}|)d^\times z.
\endaligned$$
As $S_\Lambda(y) S_\Lambda(y)=S_\Lambda(y)$, we obtain that
$$\aligned &k(x,y)=S_\Lambda(y)\{\sqrt{|{x\over y}|}h({x\over y})\\
&-\sqrt{|x|}\int_{C_S}g(|{x\over z}|)
 \left(E_S\frak F_S E_S^{-1}\right)_u\left[S_\Lambda({1\over |u|})
 \sqrt{|u|}\int_{\Bbb A_S}g(\lambda)
S_\Lambda(\lambda z)\Psi_S(\lambda zu)d\lambda\right](y)d^\times z\}.
 \endaligned $$

 By [5, Lemma 3.4 and Remark, p. 2467],
 $$\aligned &\left(e_S\frak F_S e_S^{-1}\right)_u[S_\Lambda({1\over |u|})
 \int_{\Bbb A_S}g(\lambda)
S_\Lambda(\lambda z)\Psi_S(\lambda zu)d\lambda](y)\\
&=\left(\frak F_S\right)_u[S_\Lambda({1\over |u|})
 \int_{\Bbb A_S}g(\lambda)
S_\Lambda(\lambda z)\Psi_S(\lambda zu)d\lambda](y).
\endaligned $$
It follows that
  $$\aligned &\left(E_S\frak F_S E_S^{-1}\right)_u[S_\Lambda({1\over |u|})
 \sqrt{|u|}\int_{\Bbb A_S}g(\lambda)
S_\Lambda(\lambda z)\Psi_S(\lambda zu)d\lambda](y)\\
&=\sqrt{|y|}\left(\frak F_S\right)_u[S_\Lambda({1\over |u|})
 \int_{\Bbb A_S}g(\lambda)
S_\Lambda(\lambda z)\Psi_S(\lambda zu)d\lambda](y)\\
&=\sqrt{|y|}\left(\frak F_S\right)_u[S_\Lambda({1\over |u|})
 \int_{\Bbb A_S}{1\over |z|}g({\lambda\over z})
S_\Lambda(\lambda)\Psi_S(\lambda u)d\lambda](y).\endaligned$$
Hence, by the 2nd stated identity in Lemma 4.4 we have
 $$\aligned &k(x,y)=S_\Lambda(y)\{\sqrt{|{x\over y}|}h({x\over y})\\
&-\sqrt{|xy|}\int_{C_S}g(|{x\over z}|)d^\times z
\int_{\Bbb A_S, |u|<\Lambda}\Psi_S(-uy)du
\int_{\Bbb A_S}{1\over |z|}g({\lambda\over z})
S_\Lambda(\lambda)\Psi_S(\lambda u)d\lambda\}.
 \endaligned $$
 Since
 $$\aligned \int_{C_S}g(|{x\over z}|)&d^\times z
\int_{\Bbb A_S, |u|<\Lambda}\Psi_S(-uy)du
\int_{\Bbb A_S}{1\over |z|}g({\lambda\over z})
S_\Lambda(\lambda)\Psi_S(\lambda u)d\lambda\\
&=\int_{\Bbb A_S, |u|<\Lambda}\Psi_S(-uy)du\int_{C_S}g(|{x\over z}|)d^\times z
\int_{\Bbb A_S}{1\over |z|}g({\lambda\over z})
S_\Lambda(\lambda)\Psi_S(\lambda u)d\lambda\\
&=\int_{\Bbb A_S, |u|<\Lambda}\Psi_S(-uy)du
\int_{\Bbb A_S}{1\over |x|}h({\lambda\over x})
S_\Lambda(\lambda)\Psi_S(\lambda u)d\lambda
\endaligned$$
where the changing of the order of integration is permissible as
$$\aligned \int_{C_S}&g(|{x\over z}|)d^\times z
\int_{\Bbb A_S}{1\over |z|}g({\lambda\over z})
S_\Lambda(\lambda)\Psi_S(\lambda u)d\lambda\\
&=\int_{\sqrt\mu|x|}^{|x|\over\sqrt\mu}g({|x|\over |z|}|)d^\times |z|
\int_{\Bbb A_S}{1\over |z|}g({\lambda\over |z|})
S_\Lambda(\lambda)\Psi_S(\lambda u)d\lambda.
\endaligned$$

 we have
$$k(x,x)=h(1)S_\Lambda(x)-S_\Lambda(x)|x|\int_{\Bbb A_S, |u|<\Lambda}\Psi_S(-ux)du
\int_{\Bbb A_S}{1\over |x|}h({\lambda\over x})
S_\Lambda(\lambda)\Psi_S(\lambda u)d\lambda.$$
That is,
$$\aligned
k(x,x)&=h(1)S_\Lambda(x)-S_\Lambda(x)\int_{\Bbb A_S, |u|<\Lambda |x|}\Psi_S(-u)du
\int_{\Bbb A_S}h(\lambda)S_\Lambda(\lambda x)\Psi_S(\lambda u)d\lambda\\
&=S_\Lambda(x)\int_{\Bbb A_S, |u|\geq\Lambda |x|}\Psi_S(-u)du
\int_{\Bbb A_S}h(\lambda)S_\Lambda(\lambda x)\Psi_S(\lambda u)d\lambda.
\endaligned$$

This completes the proof of the lemma.
\qed\enddemo

\demo{Proof of Theorem 1.4} By Lemma 4.5,
$$\aligned \text{trace}_{L^2(C_S)}&\left\{V(h)T\right\}
=\int_{C_S}d^\times x\int_{\Bbb A_S, |u|\geq\Lambda |x|}\Psi_S(u)du
\int_{\Bbb A_S}h(\lambda)S_\Lambda(x)S_\Lambda(\lambda x)\Psi_S(\lambda u)d\lambda\\
&=\int_{\Bbb A_S}\Psi_S(u)du\int_{C_S, |x|\leq {|u|\over\Lambda}}d^\times x
\int_{\Bbb A_S}h(\lambda)S_\Lambda(x)S_\Lambda(\lambda x)\Psi_S(\lambda u)d\lambda.
\endaligned$$
Since the above integrand is nonzero if and only if
$\max({1\over\Lambda}, {1\over\Lambda|\lambda|})<|x|\leq {|u|\over\Lambda}$,
we can write
$$\aligned &\int_{C_S, |x|\leq {|u|\over\Lambda}}d^\times x
\int_{\Bbb A_S}h(\lambda)S_\Lambda(x)S_\Lambda(\lambda x)\Psi_S(\lambda u)d\lambda\\
&=\int_{\Bbb A_S}h(\lambda)\Psi_S(\lambda u)d\lambda
\int_{C_S, \max({1\over\Lambda},{1\over\Lambda|\lambda|})<|x|\leq {|u|\over\Lambda}}
S_\Lambda(x)S_\Lambda(\lambda x)d^\times x\\
&=\int_{\Bbb A_S}h(\lambda)\Psi_S(\lambda u)
\log^+{|u|\over \max(1, {1\over|\lambda|})}d\lambda.
\endaligned$$
Hence,
$$\text{trace}_{L^2(C_S)}\left\{V(h)T\right\}
=\int_{\Bbb A_S}\Psi_S(u)du\int_{\Bbb A_S}h(\lambda)\Psi_S(\lambda u)
\log^+{|u|\over\max(1, {1\over|\lambda|})}d\lambda.$$

 Notice that
$$\aligned \int_{\Bbb A_S}&\Psi_S(v)dv\int_{\Bbb A_S}h(u)\Psi_S(u v)
\log^+{|v|\over\max(1, {1\over|u|})}du\\
=&\int_{\Bbb A_S, |v|\geq 1}\Psi_S(v)
\log|v|dv\int_{\Bbb A_S, |u|\geq 1}h(u)\Psi_S(-uv)du\\
&+\int_{\Bbb A_S}\Psi_S(v)
dv\int_{\Bbb A_S,|uv|\geq 1, |u|< 1}
h(u)\Psi_S(-uv)\log |uv|du.
\endaligned$$
We can write
$$\aligned &\int_{\Bbb A_S}\Psi_S(v)dv\int_{\Bbb A_S,|uv|\geq 1, |u|< 1}
h(u)\Psi_S(-uv)\log |uv|du\\
&=\int_{\Bbb A_S}\Psi_S(v)dv\int_{\Bbb A_S,|u|\geq 1, |u/v|< 1}
{1\over |v|}h(u/v)\Psi_S(-u)\log |u|du\\
&=\int_{\Bbb A_S}\Psi_S(vu)dv\int_{\Bbb A_S,|u|\geq 1, |1/v|< 1}
{1\over |v|}h(1/v)\Psi_S(-u)\log |u|du\\
&=\int_{|u|\geq 1}\Psi_S(-u)\log|u|du\int_{|v|>1}h(v)\Psi_S(vu)dv\\
&=\int_{|v|\geq 1}\Psi_S(v)\log|v|dv\int_{|u|>1}h(u)\Psi_S(-uv)du.
\endaligned$$
Hence,
$$\aligned \text{trace}_{L^2(C_S)}\left\{V(h)T\right\}
=&\int_{\Bbb A_S, |v|>1}\Psi_S(v)
\log|v|dv\int_{\Bbb A_S, |u|\geq 1}h(u)\Psi_S(-uv)du\\
&+\int_{|v|>1}\Psi_S(v)\log|v|dv\int_{|u|>1}h(u)\Psi_S(-uv)du.
\endaligned$$

This completes the proof of Theorem 1.4.
\qed\enddemo

\vskip 0.35in
\heading
5.  Proof of Corollary 1.5
\endheading

\proclaim{Lemma 5.1}  We can write
  $$ \aligned \int_{\Bbb A_S, |v|>1}&\Psi_S(v)
 \log|v|dv\int_{\Bbb A_S, |u|\geq 1}h(u)\Psi_S(-uvdu\\
   &=4\sum_{k_1, k_2, l_1, l_2\in\Bbb N_S}{\mu(k_1)\mu(k_2)\over k_1k_2}
  \int_1^\infty \log v\,\cos{2\pi l_1v\over k_1}dv
  \int_1^\infty h(u)\cos{2\pi l_2uv\over k_2}du.
 \endaligned $$
 \endproclaim

 \demo{Proof}  As $|\beta|=1$, $|\gamma|=1$ for $\beta,\gamma\in O_S^*$,
 $$\aligned &\int_{\Bbb A_S, |v|>1}\Psi_S(v)
\log|v|dv\int_{\Bbb A_S, |u|\geq 1}h(u)\Psi_S(-uvdu\\
&=\sum_{\beta\in O_S^*}\int_{\beta I_S, |v|>1}\Psi_S(v)\log|v|dv
\sum_{\gamma\in O_S^*}\int_{\gamma I_S, |u|\geq 1}h(u)\Psi_S(-uvdu\\
&=\sum_{\beta\in O_S^*}\varpi(\beta)\int_1^\infty e^{-2\pi iv\beta}\log|v|dv
\sum_{\gamma\in O_S^*}\varpi(\gamma)\int_1^\infty h(u)e^{2\pi\gamma uv}du
\endaligned $$
where
 $$\varpi(\beta)=\prod_{p\in S^\prime} \cases 1-p^{-1}
 &\text{if $|\beta|_p\leq 1$,}\\
 -p^{-1} &\text{if $|\beta|_p=p,$} \\
 0 &\text{if $|\beta|_p>p$}. \endcases$$
 That is,
 $$ \aligned &\int_{\Bbb A_S, |v|>1}\Psi_S(v)
\log|v|dv\int_{\Bbb A_S, |u|\geq 1}h(u)\Psi_S(-uvdu\\
 &=4\rho_S^2\sum_{k_1, l_1\in\Bbb N_S, (k_1, l_1)=1}{\mu(k_1)\over\prod_{p|k_1}(p-1)}
 \int_1^\infty \log v\,\cos{2\pi l_1v\over k_1}dv\\
 &\times\sum_{k_2, l_2\in\Bbb N_S, (k_2, l_2)=1}{\mu(k_2)\over\prod_{p|k_2}(p-1)}
  \int_1^\infty h(u)\cos{2\pi l_2uv\over k_2}du
 \endaligned \tag 5.1$$

 By partial integration,
  $$\aligned &\int_1^\infty h(u)\cos{2\pi l_2uv\over k_2}du\\
 &={k_2\over 2\pi l_2v}\left(-h(1)\sin{2\pi l_2v\over k_2}
 -\int_1^\infty h^\prime(u)\sin{2\pi l_2uv\over k_2}du\right)dv\\
 &=-h(1){k_2\over 2\pi l_2v}\sin{2\pi l_2v\over k_2}
 -h^\prime(1)({k_2\over 2\pi l_2v})^2\cos{2\pi l_2v\over k_2}
 -({k_2\over 2\pi l_2v})^2\int_1^\infty h^{\prime\prime}(u)\cos{2\pi l_2uv\over k_2}du
 \endaligned \tag 5.2$$
 and
 $$\aligned &\int_1^\infty \log v\,\cos{2\pi l_1v\over k_1}\left(
 \sum_{\underset{k_1l_2\neq k_2l_1}\to{k_2, l_2\in\Bbb N_S, (k_2, l_2)=1}}
 {\mu(k_2)\over\prod_{p|k_2}(p-1)}
 h(1){-k_2\over 2\pi l_2v}\sin{2\pi l_2v\over k_2}\right)dv\\
 &=\int_1^\infty {\log v\over v}(\sum_{\underset{k_1l_2\neq k_2l_1}
 \to{k_2, l_2\in\Bbb N_S, (k_2, l_2)=1}}
 {\mu(k_2)\over\prod_{p|k_2}(p-1)}
 h(1){-k_2\over 4\pi l_2}\left[\sin 2\pi v({l_2\over k_2}+{l_1\over k_1})
 +\sin 2\pi v({l_2\over k_2}-{l_1\over k_1})\right])dv\\
 &=\int_1^\infty {\log v-1\over v^2}\sum_{\underset{k_1l_2\neq k_2l_1}
 \to{k_2, l_2\in\Bbb N_S, (k_2, l_2)=1}}{\mu(k_2)\over\prod_{p|k_2}(p-1)}
 {h(1)k_2\over 4\pi l_2}\{{\cos 2\pi v({l_2\over k_2}-{l_1\over k_1})
 \over 2\pi({l_2\over k_2}-{l_1\over k_1})}+
 {\cos 2\pi v({l_2\over k_2}+{l_1\over k_1})
 \over 2\pi({l_2\over k_2}+{l_1\over k_1})}\}dv\\
 &=\sum_{\underset{k_1l_2\neq k_2l_1}
 \to{k_2, l_2\in\Bbb N_S, (k_2, l_2)=1}}{\mu(k_2)\over\prod_{p|k_2}(p-1)}
 {h(1)k_2\over 4\pi l_2}\int_1^\infty\{{\cos 2\pi v({l_2\over k_2}-{l_1\over k_1})
 \over 2\pi({l_2\over k_2}-{l_1\over k_1})}+
 {\cos 2\pi v({l_2\over k_2}+{l_1\over k_1})
 \over 2\pi({l_2\over k_2}+{l_1\over k_1})}\}{\log v-1\over v^2}dv\\
 &=\sum_{\underset{k_1l_2\neq k_2l_1}\to{k_2, l_2\in\Bbb N_S, (k_2, l_2)=1}}
 {\mu(k_2)\over\prod_{p|k_2}(p-1)}h(1){-k_2\over 4\pi l_2}\int_1^\infty
 {\log v\over v}\cos{2\pi l_1v\over k_1}\sin{2\pi l_2v\over k_2}dv.
 \endaligned\tag 5.3$$
 Note that
  $$\aligned &\int_1^\infty \log v\,\cos{2\pi l_1v\over k_1}
 h(1){-k_2\over 2\pi l_2v}\sin{2\pi l_2v\over k_2}dv\\
  &=\cases {h(1)k_2\over 4\pi l_2}
  \int_1^\infty\{{\cos 2\pi v({l_2\over k_2}-{l_1\over k_1})
 \over 2\pi({l_2\over k_2}-{l_1\over k_1})}+
 {\cos 2\pi v({l_2\over k_2}+{l_1\over k_1})
 \over 2\pi({l_2\over k_2}+{l_1\over k_1})}\}{\log v-1\over v^2}dv
 &\text{ if $k_2l_1\neq k_1l_2$}
 \\{h(1)k_2\over 4\pi l_2}\int_1^\infty{\cos 2\pi v({l_2\over k_2}+{l_1\over k_1})
 \over 2\pi({l_2\over k_2}+{l_1\over k_1})}{\log v-1\over v^2}dv
 &\text{ if $k_2l_1=k_1l_2$}. \endcases
 \endaligned\tag 5.4$$

  Since $\Bbb N_S$ contains only finitely many square free integers $k_1, k_2$,
 there are only finitely many integers $l_1, l_2$ in $\Bbb N_S$ such that
 $k_1l_2-k_2l_1=0$ for some square free integers $k_1, k_2$ in $\Bbb N_S$.
 It follows from (5.3) that
  $$\aligned &\int_1^\infty \log v\,\cos{2\pi l_1v\over k_1}\left(
 \sum_{k_2, l_2\in\Bbb N_S, (k_2, l_2)=1}{\mu(k_2)\over\prod_{p|k_2}(p-1)}
 h(1){-k_2\over 2\pi l_2v}\sin{2\pi l_2v\over k_2}\right)dv\\
  &=\sum_{k_2, l_2\in\Bbb N_S, (k_2, l_2)=1}
 {\mu(k_2)\over\prod_{p|k_2}(p-1)}h(1){-k_2\over 4\pi l_2}\int_1^\infty
 {\log v\over v}\cos{2\pi l_1v\over k_1}\sin{2\pi l_2v\over k_2}dv.
 \endaligned\tag 5.5$$
 From (5.2) and (5.5) we deduce that
 $$\aligned \int_1^\infty &\log v\,\cos{2\pi l_1v\over k_1}dv
 \sum_{k_2, l_2\in\Bbb N_S, (k_2, l_2)=1}{\mu(k_2)\over\prod_{p|k_2}(p-1)}
  \int_1^\infty h(u)\cos{2\pi l_2uv\over k_2}du\\
  &=\sum_{k_2, l_2\in\Bbb N_S, (k_2, l_2)=1}{\mu(k_2)\over\prod_{p|k_2}(p-1)}
  \int_1^\infty\log v\,\cos{2\pi l_1v\over k_1}dv
  \int_1^\infty h(u)\cos{2\pi l_2uv\over k_2}du.
  \endaligned \tag 5.6$$
  By (5.1) and (5.6),
  $$ \aligned \int_{\Bbb A_S, |v|>1}\Psi_S(v)
&\log|v|dv\int_{\Bbb A_S, |u|\geq 1}h(u)\Psi_S(-uvdu\\
 =&4\rho_S^2\sum_{k_1, l_1\in\Bbb N_S, (k_1, l_1)=1}
 {\mu(k_1)\over\prod_{p|k_1}(p-1)}
 \sum_{k_2, l_2\in\Bbb N_S, (k_2, l_2)=1}{\mu(k_2)\over\prod_{p|k_2}(p-1)}\\
  &\times\int_1^\infty \log v\,\cos{2\pi l_1v\over k_1}dv
  \int_1^\infty h(u)\cos{2\pi l_2uv\over k_2}du
 \endaligned \tag 5.7$$

   By [11, Example 10, p. 162],
 $$\int_0^\infty t^{s-1}\cos t\,dt=\Gamma(s)\cos{\pi s\over 2}$$
 for $0<\Re s<1$.  Since $h(u)=0$ for $u\not\in(\mu, \mu^{-1})$, we
 can change the order of integration and derive that for $0<\Re s<1$
 $$\aligned &\int_0^\infty v^{s-1}dt\int_1^\infty h(u)\cos{2\pi l_2uv\over k_2}du
 =\int_1^\infty h(u)du\int_0^\infty v^{s-1}\cos{2\pi l_2uv\over k_2}dv\\
 &=k_2^s(2l_2\pi)^{-s}\int_1^\infty h(u)u^{-s}du
 \int_0^\infty v^{s-1}\cos v\,dv\\
 &=k_2^s(2l_2\pi)^{-s}\Gamma(s)\cos{\pi s\over 2}\int_1^\infty h(u)u^{-s}du.
 \endaligned $$
 Since
 $$k_2^s(2l_2\pi)^{-s}\Gamma(s)\cos{\pi s\over 2}\int_1^\infty h(u)u^{-s}du$$
 is analytic for $0<\Re s<1$ and is $\ll |s|^{\Re s-3/2}$ in this vertical strip
  when $|s|\to\infty$. By Mellin's inversion formula [12, Theorem 28, p. 46],
  $$\aligned\int_1^\infty h(u)\cos{2\pi l_2uv\over k_2}du
  &={1\over 2\pi i}\int_{\Re s=c}k_2^s(2l_2\pi)^{-s}\Gamma(s)\cos{\pi s\over 2}
  \left(\int_1^\infty h(u)u^{-s}du\right)v^{-s}ds\\
  &={1\over 4\pi i}\int_{c-i\infty}^{c+i\infty}v^{-s}\chi(1-s)
  ({k_2\over l_2})^s\left(\int_1^\infty h(u)u^{-s}du\right)ds
  \endaligned $$
  for $0<c<1$, where $\chi(s)=2^s\pi^{s-1}\Gamma(1-s)\sin{\pi s\over 2}$.

  Similarly, by (5.2) and (5.4) we find that
  $$\aligned  \int_1^\infty &\log v\,\cos{2\pi l_1v\over k_1}dv
  \int_1^\infty h(u)\cos{2\pi l_1uv\over k_1}du
  ={1\over 4\pi i}\int_{c-i\infty}^{c+i\infty}\chi(1-z)({k_1\over l_1})^z\\
  &\times\left(\int_1^\infty \log v
  \{{1\over 4\pi i}\int_{c-i\infty}^{c+i\infty}v^{-s}\chi(1-s)
  ({k_2\over l_2})^s(\int_1^\infty h(u)u^{-s}du)ds\}v^{-z}dv\right)dz
  \endaligned \tag 5.8$$
  for $0<c<1$.

  Since $\Bbb N_S$ contains only finitely many primes and since
  $|\Gamma(x+iy)|\sim \sqrt{2\pi}|y|^{x-1/2}e^{-\pi|y|/2}$
 for fixed $x$ when $|y|\to\infty$,  the double series
  and the double integral
   $$ \aligned &4\rho_S^2\sum_{k_1, l_1\in\Bbb N_S, (k_1, l_1)=1}
 {\mu(k_1)\over\prod_{p|k_1}(p-1)}
 \sum_{k_2, l_2\in\Bbb N_S, (k_2, l_2)=1}{\mu(k_2)\over\prod_{p|k_2}(p-1)}
  {1\over 4\pi i}\int_{c-i\infty}^{c+i\infty}\chi(1-z)({k_1\over l_1})^z\\
  &\times\left(\int_1^\infty \log v
  \{{1\over 4\pi i}\int_{c-i\infty}^{c+i\infty}v^{-s}\chi(1-s)
  ({k_2\over l_2})^s(\int_1^\infty h(u)u^{-s}du)ds\}v^{-z}dv\right)dz
 \endaligned $$
 are absolute convergent for $0<c<1/2$.  Thus, by (5.7) and (5.8)
  $$ \aligned &\int_{\Bbb A_S, |v|>1}\Psi_S(v)
 \log|v|dv\int_{\Bbb A_S, |u|\geq 1}h(u)\Psi_S(-uvdu\\
 &=4\rho_S^2{1\over 4\pi i}\int_{c-i\infty}^{c+i\infty}\chi(1-z)
 [\sum_{k_1, l_1\in\Bbb N_S, (k_1, l_1)=1}
 {\mu(k_1)\over\prod_{p|k_1}(p-1)}({k_1\over l_1})^z]\int_1^\infty \log v\\
  &\times\{{1\over 4\pi i}\int_{c-i\infty}^{c+i\infty}v^{-s}\chi(1-s)
  (\sum_{k_2, l_2\in\Bbb N_S, (k_2, l_2)=1}{\mu(k_2)\over\prod_{p|k_2}(p-1)}
  ({k_2\over l_2})^s)(\int_1^\infty h(u)u^{-s}du)ds\}v^{-z}dv\}dz
 \endaligned $$
for $0<c<1/2$.

As
$$\aligned \rho_S &\sum_{k_1, l_1\in\Bbb N_S, (k_1, l_1)=1}
 {\mu(k_1)\over\prod_{p|k_1}(p-1)}({k_1\over l_1})^z
 =\sum_{k_1\in\Bbb N_S}\prod_{p|k_1}{-1\over p^{1-z}}
 \cdot \prod_{p\nmid k_1}{1-p^{-1}\over 1-p^{-z}}\\
 &=\prod_{p\in S^\prime}{1-p^{-1}\over 1-p^{-z}}
 \sum_{k_1\in\Bbb N_S}\prod_{p|k_1}{-(1-p^{-z})\over p^{1-z}(1-p^{-1})}\\
 &=\prod_{p\in S^\prime}{1-p^{-1}\over 1-p^{-z}}
 \sum_{p\in S^\prime}\left(1+{-(1-p^{-z})\over p^{1-z}(1-p^{-1})}\right)\\
 &=\prod_{p\in S^\prime}{1-p^{-1}\over 1-p^{-z}}
 \sum_{p\in S^\prime}{1-p^{z-1}\over 1-p^{-1}}
 =\prod_{p\in S^\prime}{1-p^{z-1}\over 1-p^{-z}},
 \endaligned $$
 we have
  $$ \aligned &\int_{\Bbb A_S, |v|>1}\Psi_S(v)
 \log|v|dv\int_{\Bbb A_S, |u|\geq 1}h(u)\Psi_S(-uvdu\\
 &={1\over 2\pi i}\int_{c-i\infty}^{c+i\infty}\chi(1-z)
 \prod_{p\in S^\prime}{1-p^{z-1}\over 1-p^{-z}}[\int_1^\infty \log v\\
  &\times\{{1\over 2\pi i}\int_{c-i\infty}^{c+i\infty}v^{-s}\chi(1-s)
  \prod_{p\in S^\prime}{1-p^{s-1}\over 1-p^{-s}}
  (\int_1^\infty h(u)u^{-s}du)ds\}v^{-z}dv]dz\\
  &=\sum_{k_1, k_2, l_1, l_2\in\Bbb N_S}{\mu(k_1)\mu(k_2)\over k_1k_2}
  {1\over 2\pi i}\int_{c-i\infty}^{c+i\infty}\chi(1-z)({k_1\over l_1})^z
  \int_1^\infty \log v\\
  &\times\{{1\over 2\pi i}\int_{c-i\infty}^{c+i\infty}v^{-s}\chi(1-s)
  ({k_2\over l_2})^s(\int_1^\infty h(u)u^{-s}du)ds\}v^{-z}dv\}dz\\
  &=4\sum_{k_1, k_2, l_1, l_2\in\Bbb N_S}{\mu(k_1)\mu(k_2)\over k_1k_2}
  \int_1^\infty \log v\,\cos{2\pi l_1v\over k_1}dv
  \int_1^\infty h(u)\cos{2\pi l_2uv\over k_2}du.
 \endaligned $$

 This completes the proof of the lemma.
  \qed\enddemo

\proclaim{Lemma 5.2} ([4, Theorem 8.4, p.101]) If $A$ is
of trace class on a Hilbert space, then the functional trace of $A$
coincides with its spectral trace.
\endproclaim

 \demo{Proof of Corollary 1.5} By Theorem 2.4, all eigenvalues of $V(h)T$
 on $L^2(C_S)$ are nonnegative.  Also, by Theorem 3.8 the operator
 $V(h)T$ is of trace class on the Hilbert space $L^2(C_S)$.  It follows
 from Lemma 5.2 that the functional trace of $V(h)T$ is nonnegative.

   As the functional trace of $V(h)T$ on $L^2(C_S)$ is equal to
$$\aligned \int_{\Bbb A_S, |v|>1}&\Psi_S(v)
\log|v|dv\int_{\Bbb A_S, |u|\geq 1}h(u)\Psi_S(-uvdu\\
&+\int_{|v|>1}\Psi_S(v)\log|v|dv\int_{|u|>1}h(u)\Psi_S(-u(|v|,1,\cdots,1))du
\endaligned$$
by Theorem 4.6, from Lemma 5.1 we derive that
$$\sum_{k_1, k_2, l_1, l_2\in\Bbb N_S}{\mu(k_1)\mu(k_2)\over k_1k_2}
  \int_1^\infty \log v\,\cos{2\pi l_1v\over k_1}dv
  \int_1^\infty h(u)\cos{2\pi l_2uv\over k_2}du\geq 0$$
  for all $g\in C_c^\infty(0,\infty)$.

This completes the proof of Corollary 1.5.
\qed\enddemo

\vskip 0.35in
\heading
6.  A relation to the Weil distribution
\endheading

\proclaim{Lemma 6.1} ([5, line 1-3 from bottom on p. 2483 and (3.15) on p. 2484]) We have
$$\int_{\Bbb A_S} \frak F_Sh(y)\Psi_S(y)\log|y|\,dy
=-\sum_p\int_{Q_p^*}^\prime {h(u^{-1})\over |1-u|_p}d^*u $$
where the sum on $p$ is over all primes of $Q$ including the infinity prime.
\endproclaim

We denote $x=(x_r, x_b)$ with $x_r\in\Bbb R$ and $x_b\in \Bbb A_{S^\prime}$.

\proclaim{Lemma 6.2} Let $f(\lambda)=h(\lambda )\log\max(1,{1\over|\lambda |})$.  Then
 $$\int_{\Bbb A_S}\Psi_S(-uy)du
\int_{\Bbb A_S}f(\lambda)\Psi_S(\lambda u)d\lambda=f(y).$$
\endproclaim

\demo{Proof}  By [5, Lemma 3.4 and Remark, p. 2467],
 $$\left(e_S\frak F_S e_S^{-1}\right)_u
 [\int_{\Bbb A_S}f(\lambda)\Psi_S(\lambda u)d\lambda](y)
=\left(\frak F_S\right)_u[\int_{\Bbb A_S}f(\lambda)
\Psi_S(\lambda u)d\lambda](y).$$
That is,
  $$\aligned &\left(E_S\frak F_S E_S^{-1}\right)_u
  [\sqrt{|u|}\int_{\Bbb A_S}f(\lambda)
  \Psi_S(\lambda u)d\lambda](y)\\
&=\sqrt{|y|}\left(\frak F_S\right)_u[\int_{\Bbb A_S}f(\lambda)
\Psi_S(\lambda u)d\lambda](y)\\
&=\sqrt{|y|}\int_{\Bbb A_S}\Psi_S(-uy)du\int_{\Bbb A_S}f(\lambda)
\Psi_S(\lambda u)d\lambda.
\endaligned$$

Since
$$\aligned
\int_{\Bbb A_S}f(\lambda)\Psi_S(\lambda u)d\lambda
&=\sum_{\gamma\in O_S^*}\int_{\gamma I_S}
f(\lambda)\Psi_S(\lambda u)d\lambda\\
&=\sum_{\gamma\in O_S^*}\int_{I_S}f(\lambda)
\Psi_S(\lambda\gamma u)d\lambda,
\endaligned $$
we have
$$E_S^{-1}\{\sqrt{|u|}\int_{\Bbb A_S}f(\lambda)
\Psi_S(\lambda u)d\lambda\}
=\int_{I_S}f(\lambda)\Psi_S(\lambda u)d\lambda.$$
Hence,
$$\aligned
&(\frak F_SE_S^{-1})_u\{\sqrt{|u|}
\int_{\Bbb A_S}f(\lambda)\Psi_S(\lambda u)d\lambda\}(y)\\
&=\int_{\Bbb A_S}\Psi_S(-uy)du
\int_{I_S}f(\lambda)\Psi_S(\lambda u)d\lambda\\
&=\int_{-\infty}^\infty e^{2\pi iu_ry_r}du_r \int_{\Bbb A_{S^\prime}}
\Psi_{S^\prime}(-u_by_b)du_b\int_{\Bbb A_{S^\prime}}\phi(\lambda_b)
\Psi_{S^\prime}(\lambda_bu_b)d\lambda_b \endaligned $$
where
$$\phi(\lambda_b)=\cases \int_{-\infty}^\infty f(|\lambda_r||\lambda_b|)
e^{-2\pi i\lambda_ru_r}d\lambda_r
 &\text{if $\lambda_b\in\prod_{p\in S^\prime}O_p^*$}\\
0&\text{if $\lambda_b\not\in\prod_{p\in S^\prime}O_p^*$.}\endcases$$
Also,
$$\aligned &\int_{I_S}f(\lambda)\Psi_S(\lambda u)d\lambda
=\varpi(u)\int_0^\infty f(t)e^{-2\pi itu_r}dt.
\endaligned $$
This implies that
$$\int_{\Bbb A_{S^\prime}}\phi(\lambda_b)
\Psi_{S^\prime}(\lambda_bu_b)d\lambda_b$$
as a function of $u_b$ is supported on the compact set
$B=\prod_{p\in S^\prime}p^{-1}O_p$ with $O_p=\{x\in Q_p: |x|_p\leq 1\}$.
Since $\phi(\lambda_b)$ is locally constant as
a function of $\lambda_b\in\Bbb A_{S^\prime}$ and is supported on the compact set
$\prod_{p\in S^\prime}O_p^*$,
the condition of [10, Theorem 2.2.2, p. 310] is satisfied by this function.
By the Fourier inversion formula (Note that [10, Theorem 2.2.2, p. 310] is still
true if we replace $k_p^+$ there by $\Bbb A_{S^\prime}$),  we have
$$\int_{\Bbb A_{S^\prime}}
\Psi_{S^\prime}(-u_by_b)du_b\int_{\Bbb A_{S^\prime}}\phi(\lambda_b)
\Psi_{S^\prime}(\lambda_bu_b)d\lambda_b=\phi(y_b)$$
where $\phi(y_b)=0$ if $y_b\not\in\prod_{p\in S^\prime}O_p^*$.

 Since $f(t|y_b|)$  is a continuous and compactly supported
function of $t\in(0,\infty)$, $t\neq 1/\Lambda|y_b|$
 and is of bounded variation in an interval including
$y_r$, by Fourier's single-integral formula [12, Theorem 12, p. 25] we have
for $y_b\in\prod_{p\in S^\prime}O_p^*$
$$\aligned
&(\frak F_SE_S^{-1})_u\{\sqrt{|u|}
\int_{\Bbb A_S}f(\lambda)\Psi_S(\lambda u)d\lambda\}(y)\\
&=\int_{-\infty}^\infty e^{2\pi iu_ry_r}du_r
\int_{-\infty}^\infty f(|\lambda_r||y_b|)
 e^{-2\pi i\lambda_ru_r}d\lambda_r\\
  &=\lim_{X\to\infty}\int_{-\infty}^\infty f(|\lambda_r||y_b|)
   d\lambda_r\int_{-X}^X e^{2\pi iu_r(y_r-\lambda_r)}du_r\\
 &=\lim_{X\to\infty}{1\over\pi}\int_{-\infty}^\infty f(|\lambda_r||y_b|)
   {\sin 2\pi X(y_r-\lambda_r)\over y_r-\lambda_r}d\lambda_r
  =f(y)
  \endaligned $$
  where the 2nd equality holds because $f$ has a compact support on $(0,\infty)$.
  It follows that
 $$\aligned
&(E_S\frak F_SE_S^{-1})_u\{\sqrt{|u|}\int_{\Bbb A_S}f(\lambda)
\Psi_S(\lambda u)d\lambda\}(y)\}(y)\\
  &=\sqrt{|y|}\sum_{\xi\in O_S^*}\cases f(\xi y)
   &\text{if $\xi y\in I_S$}, \\ 0 &\text{if $\xi y\not\in I_S$}\endcases
   =\sqrt{|y|}f(y)
\endaligned$$
because for each $y\in\Bbb A_S$ with $|y|\neq 0$ there exists exactly one
$\xi\in O_S^*$ such that $\xi y\in I_S$.  Therefore,
$$\left(E_S\frak F_S E_S^{-1}\right)_u
  [\sqrt{|u|}\int_{\Bbb A_S}f(\lambda)
\Psi_S(\lambda zu)d\lambda](y)=\sqrt{|y|}f(y).$$

This completes the proof of the lemma.
\qed\enddemo

\proclaim{Corollary 6.3}  We have
$$\int_{\Bbb A_S}\Psi_S(-v)dv
\int_{\Bbb A_S}h(y)\log\max(1,{1\over|y|})\Psi_S(yv)dy=0.$$
\endproclaim

\demo{Proof}  By Lemma 6.2,
$$\int_{\Bbb A_S}\Psi_S(-vx)dv
\int_{\Bbb A_S}h(y)\log\max(1,{1\over|y|})\Psi_S(yv)dy
=h(x)\log\max(1,{1\over|x|}).$$
Choosing $x=1$ in the above identity gives the stated identity.

This completes the proof of the lemma.
\qed\enddemo

  We denote $\log^+x=\log x$ for $x>1$, $\log^+x=0$ for $0<x\leq 1$ and
  $\log^-x=\log x$ if $0<x\leq 1$, $\log^-x=0$ if $x> 1$.

\demo{Proof of Theorem 1.6} Since
$$\int_{\Bbb A_S}\Psi_S(-v)dv
\int_{\Bbb A_S}h(u)\log\max(1,{1\over|u|})\Psi_S(uv)du=0,$$
 by Lemma 6.1 and Corollary 6.3 we can write
 $$\aligned
 &\Delta(h)=\int_{\Bbb A_S}\Psi_S(-v)dv\int_{\Bbb A_S}h(u)\log|v|\Psi_S(uv)du\\
&-\int_{\Bbb A_S}\Psi_S(-v)dv\int_{\Bbb A_S}h(u)\log\max(1,{1\over|u|})
\Psi_S(uv)du\\
&=\int_{\Bbb A_S}\Psi_S(-v)dv\int_{\Bbb A_S}h(u)\log{|v|\over\max(1,{1\over|u|})}
\Psi_S(uv)du\\
&=\int_{\Bbb A_S}\Psi_S(-v)dv\int_{\Bbb A_S}h(u)\left(\log^+{|v|\over\max(1,{1\over|u|})}
+\log^-{|v|\over\max(1,{1\over|u|})}\right)\Psi_S(uv)du\\
&=\int_{|v|>1}\log|v|\Psi_S(-v)dv\int_{|u|\geq 1}h(u)\Psi_S(uv)du
+\int_{\Bbb A_S}\Psi_S(-v)dv\int_{|vu|>1, |u|<1}h(u)\log|vu|\Psi_S(uv)du\\
&+\int_{|v|<1}\log|v|\Psi_S(-v)dv\int_{|u|\geq 1}h(u)\Psi_S(uv)du
+\int_{\Bbb A_S}\Psi_S(-v)dv\int_{|vu|<1, |u|<1}h(u)\log|vu|\Psi_S(uv)du.
 \endaligned $$
 As
 $$\aligned &\int_{\Bbb A_S}\Psi_S(-v)dv\int_{|vu|>1, |u|<1}h(u)\log|vu|\Psi_S(uv)du\\
 &=\int_{\Bbb A_S}\Psi_S(-v){dv\over |v|}\int_{|u|>1, |u|<|v|}h({u\over v})\log|u|\Psi_S(u)du\\
 &=\int_{|u|>1}\log|u|\Psi_S(u)du\int_{|u|<|v|}{1\over |v|}h({u\over v})\Psi_S(-v)dv\\
 &=\int_{|u|>1}\log|u|\Psi_S(u)du\int_{|v|>1}h(v)\Psi_S(-vu)dv
 \endaligned $$
 where we used the fact that $|v|^{-1}h(v^{-1})=h(v)$.  Similarly we have
 $$\aligned &\int_{\Bbb A_S}\Psi_S(-v)dv\int_{|vu|<1, |u|<1}h(u)\log|vu|\Psi_S(uv)du\\
&=\int_{|u|<1}\log|u|\Psi_S(u)du\int_{|v|>1}h(v)\Psi_S(-vu)dv.
\endaligned$$
Therefore,
 $$\aligned
 &\Delta(h)=\int_{|v|>1}\log|v|\Psi_S(-v)dv\int_{|u|\geq 1}h(u)\Psi_S(uv)du
+\int_{|v|>1}\log|v|\Psi_S(-v)dv\int_{|u|>1}h(u)\Psi_S(uv)du\\
&+\int_{|v|<1}\log|v|\Psi_S(-v)dv\int_{|u|\geq 1}h(u)\Psi_S(uv)du
+\int_{|v|<1}\log|v|\Psi_S(-v)dv\int_{|u|>1}h(u)\Psi_S(uv)du.
 \endaligned $$
 By Theorem 1.4,
$$\aligned \text{trace}_{L^2(C_S)}\left\{V(h)T\right\}
=&\int_{\Bbb A_S, |v|>1}\Psi_S(v)
\log|v|dv\int_{\Bbb A_S, |u|\geq 1}h(u)\Psi_S(-uv)du\\
&+\int_{|v|>1}\Psi_S(v)\log|v|dv\int_{|u|>1}h(u)\Psi_S(-uv)du.
\endaligned$$
The stated identity follows from the above two identities.

This completes the proof of Theorem 1.6.
\qed\enddemo

\enddocument